\documentclass[11pt, reqno]{amsart}

\evensidemargin
\oddsidemargin 
\makeatletter\@addtoreset{equation}{section}\makeatother

\usepackage{amsmath, amsthm, latexsym, amssymb, graphics, bm, color}

\newcommand{\eps}{\varepsilon}

\newcommand{\e}{\varepsilon}

\renewcommand{\phi}{\varphi}

\newcommand{\ssup}[1] {{\scriptscriptstyle{({#1}})}}

\newcommand{\one}{\1}

\newcommand{\Prob}{{\mathrm{Prob}}}

\newcommand{\Deltad}{\Delta}
\newcommand{\dd}{\,\mathrm{d}}
\newcommand{\ee}{\mathrm{e}}
\newcommand{\vt}{v_t}
\newcommand{\lt}{\lambda_t}

\newcommand{\R}{\mathbb R}
\newcommand{\Z}{\mathbb Z}

\newcommand{\N}{\mathbb N}

\newcommand{\E}{\mathbb E}
\renewcommand{\P}{\mathbb P}

\newtheorem{theorem}{Theorem}

\newtheorem{lemma}{Lemma}

\newtheorem{prop}{Proposition}

\newcommand{\sfrac}[2] {\mbox{$\frac{#1}{#2}$}}
\newenvironment{Proof}[1]
{\vskip0.1cm\noindent{\bf #1}{\hspace*{0.3cm}}}{\vspace{0.15cm}}

\def\1{{\mathchoice {1\mskip-4mu\mathrm l}      % Blackboard bold 1
{1\mskip-4mu\mathrm l}
{1\mskip-4.5mu\mathrm l} {1\mskip-5mu\mathrm l}}}

\renewcommand{\subsection}{\secdef \subsct\sbsect}
\newcommand{\subsct}[2][default]{\refstepcounter{subsection}
\vspace{0.15cm}
{\flushleft\bf \arabic{section}.\arabic{subsection}~\bf #1  }
\nopagebreak\nopagebreak}
\newcommand{\sbsect}[1]{\vspace{0.1cm}\noindent
{\bf #1}\vspace{0.1cm}}

\newcounter{remnr}

\newenvironment{remark}{\refstepcounter{remnr}
{\sf Remark~\arabic{remnr}.\ }\nopagebreak  }%
{\nopagebreak {\hfill{$\diamond$}}\\ }

\begin{document}
\title[Complete localisation in the parabolic Anderson model]
{\large Complete localisation in the parabolic Anderson model with Pareto-distributed potential}

\author[Wolfgang K\"onig, Peter M\"orters, and Nadia Sidorova]{}

\maketitle

\centerline{\sc Wolfgang K\"onig
\qquad\sc Peter M\"orters\qquad 
\sc Nadia Sidorova}

\vspace{0.4cm}

% \centerline{\small(\version)}

\vspace{0.4cm}

\begin{quote}
{\small {\bf Abstract:} The parabolic Anderson problem is the Cauchy problem for the heat equation 
$\partial_t u(t,z)=\Delta u(t,z)+\xi(z) u(t,z)$ on $(0,\infty)\times {\mathbb Z}^d$ with random potential 
$(\xi(z) \colon z\in {\mathbb Z}^d)$. We consider independent and identically distributed potential variables,
such that Prob$(\xi(z)>x)$ decays polynomially
as $x\uparrow\infty$. If $u$ is initially localised in the origin, i.e. if $u(0,x)=\one_0(x)$,
we show that, at any large time $t$, the solution is completely localised in a single point with high probability.
More precisely, we find a random process $(Z_t \colon t\ge 0)$ with values in $\Z^d$ such that 
$\lim_{t \uparrow\infty} u(t,Z_t)/\sum_{z\in\Z^d} u(t,z) =1,$ in probability. 
We also identify the asymptotic behaviour of $Z_t$ in terms of a weak limit theorem.}
\end{quote}

\renewcommand{\thefootnote}{}
\renewcommand{\thefootnote}{1}

\section{Introduction and main results}

\subsection{The parabolic Anderson model and intermittency}\label{sec-intro}\\[-3mm]

\noindent We consider the heat equation with random potential on the integer lattice $\Z^d$ and study 
the Cauchy problem with localised initial datum,
\begin{align}\label{pam}
\begin{array}{rcll}
\displaystyle \vspace{2mm} \partial_t u(t,z) & = & \Deltad u(t,z)+\xi(z)u(t,z), \qquad
& (t,z)\in (0,\infty)\times \Z^d,\\
u(0,z) & = & \one_{0}(z), & z\in\Z^d,
\end{array}
\end{align}
where 
\begin{align*}
(\Deltad f)(z)=\sum_{y\sim z}[f(y)-f(z)], \qquad z\in\Z^d,\, f\colon \Z^d\to\R,
\end{align*}
is the discrete Laplacian, and the potential 
$(\xi(z) \colon z\in\Z^d)$ is a collection 
of independent identically distributed random variables. 
\medskip

The problem~\eqref{pam} and its variants are often called the 
\emph{parabolic Anderson problem}. The elliptic version of this problem
originated in the work of the physicist P.~W.~Anderson on entrapment of electrons
in crystals with impurities, see~\cite{An58}. The parabolic version of the problem appears
in the context of chemical kinetics and population dynamics, and also provides a simplified 
qualitative approach to problems in magnetism and turbulence.
The references \cite{GM90}, \cite{M94} and \cite{CM94} provide 
applications, background and heuristics around the parabolic Anderson model.
Interesting recent mathematical progress can be found, for example in~\cite{BMR05}, 
\cite{HKM06}, and \cite{GH06}, and \cite{GK05} is a recent survey article.\medskip

% besser [HKM06] durch neues Ergebnis anderer Autoren ersetzen,
% besser [BMR05] durch neues Cranston Ergebnis ersetzen, wenn verfuegbar

One main reason for the great interest in the parabolic Anderson problem lies in the fact that it 
exhibits an \emph{intermittency effect}: It is believed that, at late times,  the overwhelming contribution to 
the total mass of the solution $u$ of the problem~\eqref{pam} comes from a small number of
widely separated regions of small diameter, which are often called the {\it relevant islands}.
This effect is believed to get stronger (with a smaller number of relevant islands, which are of smaller size) as
the tail of the potential variable at infinity gets heavier. Providing rigorous evidence for 
intermittency is a major challenge for mathematicians, which
has lead to substantial research efforts in the past 15 years.\medskip

An approach, which has been proposed in the physics literature, see~\cite{ZM+87} or \cite{GK05},
suggests to study large time asymptotics of the moments of the total mass
\begin{equation}\label{U(t)def}
U(t)=\sum_{z\in\Z^d}u(t,z) , \qquad t>0\, .
\end{equation}
Denoting expectation with respect to $\xi$ by $\langle\,\cdot\,\rangle$, 
if all exponential moments $\langle \exp(\lambda \xi(z)) \rangle$ for $\lambda>0$ exist, then so do
all moments $\langle U(t)^p\rangle$ for $t>0$, $p>0$. Intermittency becomes manifest in a faster 
growth rate of higher moments. More precisely, the model is called intermittent~if
\begin{equation}\label{Intermitt}
\limsup_{t\to\infty}\frac {\langle U(t)^p\rangle^{1/p}}{\langle U(t)^q\rangle^{1/q}}=0,
\qquad \mbox{ for $0<p<q$. }
\end{equation}
Whenever $\xi$ is nondegenerate random, the parabolic Anderson model is intermittent in this sense, see \cite[Theorem~3.2]{GM90}. 
Further properties of the relevant islands, like their asymptotic size and
shape of potential and solution, are reflected (on a heuristical level) 
in the asymptotic expansion of $\log \langle U(t)^p\rangle$
for large $t$. Recently, in \cite{HKM06}, it was argued that the distributions with finite exponential
moments can be divided into exactly four different universality classes, with each class having a qualitatively
different long-time behaviour of the solution.\medskip

It is, however, a much harder mathematical challenge to prove intermittency in the original
geometric sense, and to identify asymptotically the number, size and location of the relevant
islands. This programme was initiated by Sznitman for the closely related continuous
model of a Brownian motion with Poissonian obstacles, and the very substantial body of research he
and his collaborators created is surveyed in his monograph~\cite{S98}. For the problem~\eqref{pam} and
two universality classes of potential distributions, the double-exponential distribution and distributions
with tails heavier than double-exponential (but still with all exponential moments finite), 
the recent paper~\cite{GKM06} makes substantial progress towards completing the geometric picture:
Almost surely, the contribution coming from the complement of a random 
number of relevant islands is negligible compared to the mass coming from these islands, asymptotically
as $t\to\infty$. In the double-exponential case, the radius of the islands stays bounded, in the heavier
case the islands are single sites, and in Sznitman's case the radius tends to infinity. 
%on the scale~$t^{1/(d+2)}$. % P removed, as I could not confirm this from PTRF 105,31-56 (1996).
\medskip

Questions about the number of relevant islands remained open in all these cases. 
Both in \cite{GKM06} and \cite{S98} it is shown that an upper bound on the number of relevant
islands is~$t^{o(1)}$, but this is certainly not always best possible. In particular, the questions 
whether a \emph{bounded number} of islands already carry the bulk of the mass, or when \emph{just one}
island is sufficient, are unanswered. These questions are difficult, since there are 
many local regions that are good candidates for being a relevant island, and the known 
criteria that identify relevant islands do not seem to be optimal.
\medskip

In the present paper, we study the parabolic Anderson model with potential distributions 
that do not have any finite exponential moment. 
For such distributions one expects the intermittency effect to be even more pronounced than in the cases 
discussed above, with a very small number of relevant islands, which are just single sites. Note that
in this case intermittency cannot be studied in terms of the moments $\langle U(t)^p \rangle$, which are not finite. 
\medskip

The main result of this paper is that, in the case of Pareto-distributed potential variables, there is only a \emph{single}
relevant island, which consists of a single site. In other words, at any large time~$t$, with high probability,
the total mass $U(t)$ is concentrated in a single lattice point $Z_t\in\Z^d$. This extreme form of intermittency
is called \emph{complete localisation}. It has been observed so far only for quite simple mean field models, see  \cite{FM, FG}, 
and the present paper is the first instance where it has been found in the parabolic Anderson model
or, indeed, any comparable lattice-based model. 
We also study the asymptotics of the location $Z_t$ of the point where the mass concentrates:  
We show that $Z_t$ goes to infinity like $(t/\log t)^{\alpha/(\alpha-d)}$, where $\alpha>d$ is
the parameter of the Pareto distribution. The location of the relevant island is further described in terms of a 
weak limit theorem for the scaled quantity $(t/\log t)^{\alpha/(d-\alpha)}Z_t$ with an explicit 
limiting density.  Precise statements are formulated in the next section.
\bigskip

\subsection{The parabolic Anderson model with Pareto-distributed potential}\label{sec-result}\\[-3mm]

We assume that the potential variables~$\xi(z)$ at all sites $z$ are independently \emph{Pareto-distributed}
with parameter $\alpha>d$, i.e., the distribution function is
\begin{equation}
F(x)=\Prob(\xi(z)\le x)=1-x^{-\alpha},\qquad x\ge 1.
\end{equation}
In particular, we have $\xi(z)\geq 1$ for all $z\in\Z^d$, almost surely. Note from \cite[Theorem~2.1]{GM90} that the 
restriction to parameters $\alpha>d$ is necessary and sufficient for \eqref{pam} to possess a 
unique nonnegative solution $u\colon (0,\infty)\times\Z^d\to[0,\infty)$. Recall that
$U(t)=\sum_{z\in\Z^d}u(t,z)$ is the total mass of the solution at time~$t>0$. We introduce
\begin{equation}\label{mudef}
q=\frac{d}{\alpha-d}\qquad\mbox{and}\qquad\mu=\frac{(\alpha-d)^d 2^d B(\alpha-d,d)}{d^d(d-1)!},
\end{equation}
where $B(\cdot,\cdot)$ denotes the Beta function. Throughout the paper we use $|x|$ 
to denote the $\ell^1$-norm of a vector $x\in\R^d$.
\medskip

Our first main result shows the complete localisation of the solution $u(t, \,\cdot\,)$ in a single
lattice point~$Z_t$, as $t\to\infty$.

\begin{theorem}[Concentration in one site]
\label{main}
There exists a process $(Z_t \colon t>0)$ with values in~$\Z^d$ such that  
$$
\lim_{t\to\infty}\frac{u(t,Z_t)}{U(t)}= 1\qquad\mbox{in probability.}
$$
\end{theorem}
\medskip

\begin{remark}
Our statement is formulated in terms of a convergence \emph{in probability}. The
convergence does \emph{not hold} in the almost sure sense. Indeed, suppose that $t>0$ 
is sufficiently large to ensure $u(t,Z_t)\ge \frac23 U(t)$ and that $t$ is 
a jumping time, i.e. that $Z_{t-}\not= Z_{t+}$. Then, by continuity, we have
$u(t,Z_{t-}) + u(t,Z_{t+}) \ge \frac43 U(t)$, which contradicts the nonnegativity
of the solution. 
\end{remark}

\begin{remark}
The asymptotic behaviour of $\log U(t)$ for the Anderson model with heavy-tailed potential
variables is analysed in detail in \cite{HMS06}. In the case of a Pareto-distributed potential
it turns out that already the leading term in the asymptotic expansion of $\log U(t)$ is random. 
This is in sharp contract to potentials with exponential moments, where the leading two
terms in the expansion are always deterministic.
More precisely, in \cite[Theorem~1.2]{HMS06} the following limit law for $\log U(t)$ is proved, 
\begin{equation}\label{U(t)asy}
\frac{(\log t)^{q}}{t^{q+1}} \, \log U(t)\Longrightarrow Y,\qquad\mbox{ where }\quad 
\P(Y\leq y)=\exp\{-\mu y^{d-\alpha}\} \mbox{ for } y>0.
\end{equation}
Note that the upper tails of $Y$ have the same asymptotic order as the Pareto distribution with parameter $\alpha-d$, 
i.e., $\P(Y>y)\asymp y^{d-\alpha}$ as $y\to\infty$. A careful inspection of the proof of~\cite[Theorem~1.2]{HMS06} shows that also
\begin{equation}\label{u(t,Zt)law}
\frac{(\log t)^{q}}{t^{q+1}} \, \log u(t,Z_t) \Longrightarrow Y.
\end{equation}
Note, however, that a combination of \eqref{U(t)asy} with \eqref{u(t,Zt)law} does not yield the 
concentration property in Theorem~\ref{main}. Much more precise techniques are necessary.
\end{remark}

Our second main result is a limit law for the concentration site $Z_t$ in Theorem~\ref{main}. 
Recall the definition of $q$ and $\mu$ from~\eqref{mudef}.
As usual, we denote weak convergence by~$\Rightarrow$.

\begin{theorem}[Limit law for the concentration site]\label{thm-Ztlim}
As $t\to\infty$, 
\begin{equation}\label{conv}
Z_t\, \left(\frac{\log t}{t}\right)^{q+1}
\Longrightarrow X,
\end{equation}
where 
$X$ is an  $\R^d$-valued random variable with Lebesgue density
\begin{align*}
p(x)=\alpha\int_0^{\infty}\frac{\exp\{-\mu y^{d-\alpha}\}}{(y+q|x|)^{\alpha+1}}\, \dd y\, .
\end{align*}
\end{theorem} 
\smallskip
% The proof of Theorem~\ref{thm-Ztlim} is at the end of Section~\ref{sec-Ztdef}.

\begin{remark} Note that $X$ is isotropic in the
$\ell^1$-norm. \end{remark}

\begin{remark}\label{3}%
The density $p$ is a probability density. Indeed, 
$$
\begin{aligned}
\int_{\R^d} p(x)\dd x&=\alpha \int_0^{\infty}\dd y\,e^{-\mu y^{d-\alpha}}\int_{\R^d}\dd x\, {(y+q|x|)^{-(\alpha+1)}},
\end{aligned}$$
and, by \cite[Lemma~3.9]{HMS06}, the inner integral equals $\frac{2^dq^{-d}}{(d-1)!}  \, B(\alpha+1-d,d) \, y^{-\alpha+d-1}$. 
Using a change of variable and the definition of $\mu$ in \eqref{mudef}, this
simplifies to
$$\int_{\R^d} p(x)\dd x = \frac{B(\alpha+1-d,d)}{B(\alpha-d,d)}\frac{\alpha}{\alpha-d}\,\mu\int_0^\infty e^{-\mu t}\dd t.$$
The integral equals $1/\mu$, and the remaining product equals one because of the 
functional equation $(x+y) \, B(x+1,y)=x\,B(x,y)$ 
for $x,y>1$, which is satisfied by the Beta function.
Moreover, the proof of Theorem~\ref{thm-Ztlim} shows that the two limit laws in \eqref{u(t,Zt)law} and 
in \eqref{conv} hold jointly, and the joint density of $(X,Y)$ is the map 
$$(x,y)\mapsto \alpha\frac{\exp\left\{-\mu y^{d-\alpha}\right\}}{(y+q|x|)^{\alpha+1}}.$$
This explains the structure of the density $p(x)$.
\end{remark}

\subsection{Overview: The strategy behind the proofs}\label{overview} \\[-2mm]

As shown in \cite[Theorem 2.1]{GM90}, under the assumption $\alpha>d$, the 
unique nonnegative solution $u\colon (0,\infty)\times\Z^d\to[0,\infty)$
of \eqref{pam} has a \emph{Feynman-Kac representation} 
\begin{align*}
u(t,z)=\E_0 \Big[\one{\{X_t=z\}}\,  \exp\Big\{\int_0^t\xi(X_s)\dd s\Big\}\Big], \qquad t>0, \, z\in\Z^d,
\end{align*} 
where $(X_s \colon s\ge 0)$ under $\P_0$ (with expectation $\E_0$) is a
continuous-time simple random walk on $\Z^d$ with generator~$\Deltad$ started in the origin.
Hence, the total mass of the solution is given by
\begin{align*}
U(t)=\E_0 \Big[\exp\Big\{\int_0^t\xi(X_s)\dd s\Big\}\Big].
\end{align*} 
Heuristically, for a fixed, large time~$t>0$, the walks $(X_s \colon 0 \le s \le t)$ that have the greatest 
impact on the average $U(t)$ move quickly to a remote site $z$ which,
\begin{itemize}
\item has a large potential value~$\xi(z)$,
\item and can be reached quickly, i.e. is sufficiently close to the origin.
\end{itemize}
Once this site is reached, the walk remains there until time~$t$. As the probability of
moving to a site $z$ within $t$ time units is approximately
$$\P_0( X_t=z ) = \exp\big\{ - |z| \, \log\big( \sfrac{|z|}{2det}\big) \, (1+ o(1))\, \big\},$$
it is plausible that the optimal site~$z$ at time~$t$ is the maximiser~$Z_t$ of the 
random functional
$$\Psi_t(z)=\xi(z)-\frac{|z|}{t}\log\frac{|z|}{2det}\, ,$$
with the understanding that $\Psi_t(0)=\xi(0)$. 
This is indeed the definition of the process $(Z_t \colon t\ge 0)$, which is underlying our two main theorems.
\medskip

In Section~\ref{sec-Ztdef} we study the asymptotic behaviour of $(Z_t \colon t\ge 0)$
using techniques from extremal value theory. We prove Theorem~\ref{thm-Ztlim} and
also provide auxiliary results that compare the largest and second-largest value in the set
$\{ \Psi_t(z) \colon z\in\Z^d\}$, as needed in the proof of Theorem~\ref{main} 
in Section~\ref{sec-proofmain}. Note that the arguments in this section are based entirely on the definition
of $(Z_t \colon t\ge 0)$ in terms of~$\Psi_t$, and not on its r\^ole in the parabolic
Anderson problem.
\medskip

Section~\ref{sec-proofmain} is devoted to the proof of Theorem~\ref{main}. In this proof we 
build on techniques developed in~\cite{GKM06}. We split~$u$ into three terms, which correspond 
to the contributions to the Feynman-Kac formula coming from paths that (1)~by time~$t$ have left 
a centred box with a certain large, $t$-dependent, random radius, (2)~stay inside 
this box for $t$ time units but do not visit~$Z_t$, and (3)~stay inside this box and do visit~$Z_t$. 
It will turn out that the total mass of the first two contributions is negligible,
%see Proposition~\ref{lu1u2}, 
and  that the total mass of the last one is concentrated on~$Z_t$.
%see Proposition~\ref{lu3}. 
To be more precise, we denote the three parts in the decomposition 
by $$u(t,z)=u_1^{\ssup t}(t,z)+u_2^{\ssup t}(t,z)+u_3^{\ssup t}(t,z)\, .$$
The radius of the box will be chosen large enough that $u_1^{\ssup t}$ has small total 
mass relative to $U(t)$, since it is expensive to reach the complement of the large box. 
\medskip

In order to deal with $u_2^{\ssup t}$, we use the gap between the value of $\Psi_t$ in its maximum~$Z_t$, 
and the maximum of $\Psi_t(z)$ over all other points $z\in \Z^d\setminus \{Z_t\}$, i.e. the auxiliary 
result provided in Section~\ref{sec-Ztdef}. From this we infer that the total mass of $u_2^{\ssup t}$ 
is small, as the site~$Z_t$, which maximises $\Psi_t$, is ruled out from the exponential.
\medskip

Finally, for the estimate of $u_3^{\ssup t}$ it is crucial that the radius of the box 
is chosen in such a way that $Z_t$ is also a maximiser of the field $\xi$ over the box. 
The main ingredient is a spectral analytical device, which is used in a similar manner 
as in~\cite{GKM06}: We show that $u_3^{\ssup t}$ can be controlled in terms of the principal eigenfunction 
of the Anderson Hamiltonian, $\Delta+\xi$, in the box with zero boundary conditions. This eigenfunction 
turns out to be exponentially concentrated in the maximal potential point in the box, which is $Z_t$. 
Hence the total mass of~$u$ is concentrated in~$Z_t$. This argument is the key
step in the proof of Theorem~\ref{main}. %, see Lemma~\ref{lem-lu3}. 
\bigskip

\section{Proof of Theorem~\ref{thm-Ztlim}: The concentration site $Z_t$}\label{sec-Ztdef}

In this section, we study the top two values in the order statistics 
of the random variables $( \Psi_t(z) \colon z\in \Z^d )$. We first prove that,
for any $t>0$,   the set $\{ \Psi_t(z) \colon z\in \Z^d \}$ is almost surely bounded. 
Thus, by continuity of  the distribution function~$F$, the set $\{ \Psi_t(z) \colon z\in \Z^d \}$ has a
unique maximum $Z_t=Z_t^{\ssup 1}\in\Z^d$. Moreover, the set
$\{ \Psi_t(z) \colon z\in \Z^d\setminus\{Z_t^{\ssup 1}\} \}$ also has a unique maximum,
which we denote by $Z_t^{\ssup 2}\in\Z^d$. Note that $Z_t^{\ssup 1}>Z_t^{\ssup 2}$.
\medskip

\begin{lemma} For any $t>0$, $\Psi_t$ is bounded almost surely.
\end{lemma}

\begin{Proof}{Proof.}
For any $r>0$, let $\xi_r^{\ssup 1}=\max_{z\in\Z^d\colon |z|\leq r}\xi(z)$ denote the maximum of the 
potential in the box with radius $r$. Denote $\phi(x)=-\log(1-F(x))$. By~\cite[Lemma 4.2]{GM90}, 
which holds for any distribution function, we have, almost surely,
\begin{align*}
\phi(\xi_r^{\ssup 1})= d\log r(1+o(1)),\qquad\text{as }r\to\infty. 
\end{align*}
Fix $\e\in(0,1-\frac d\alpha)$. Since $\phi(x)=\alpha\log x$, there 
exists a random radius $\rho_1>0$ such that, almost surely,
\begin{align}
\label{asas}
\xi_r^{\ssup 1}\le r^{\frac{d}{\alpha}+\e},\qquad \mbox{for all }r>\rho_1. 
\end{align}

Now fix $t>0$.
Since $\frac{d}{\alpha}+\e<1$, 
there exists $\rho_{2}(t)>0$ such that
$r^{\frac{d}{\alpha}+\e}<\frac{r}{t}\log\frac{r}{2det}$ for all $r>\rho_{2}(t)$. With $\rho(t)=
\max\{\rho_1,\rho_2(t)\}$, we obtain
\begin{align*}
\sup_{|z|>\rho(t)}\Psi_t(z)
%\le \sup_z \left[\xi(z)-\frac{|z|}{t}\log\frac{|z|}{t}\right]
&\le\sup_{|z|>\rho(t)}\left[\xi_{|z|}^{\ssup 1}-\frac{|z|}{t}\log\frac{|z|}{2det}\right]
\le\sup_{r>\rho(t)}\left[r^{\frac{d}{\alpha}+\e}-\frac{r}{t}\log\frac{r}{2det}\right]\le 0.
\end{align*}
Hence, the function $\Psi_t$ is positive only for finitely many $z$ and thus attains 
its maximum.
\qed
\end{Proof}

We define two scaling functions 
\begin{equation}\label{rtatdef}
r_t=\left(\frac{t}{\log t}\right)^{q+1}
\qquad\text{and}\qquad
{a}_t=\left(\frac{t}{\log t}\right)^{q}.
\end{equation}
A limit law for $\Psi_t(Z_t)$ is given in \cite[Prop.~3.8]{HMS06}. From its proof it follows that
\begin{align}
\label{psiz}
\Prob\left(\Psi_t(Z_t)\le {a}_ty\right)
=e^{-\mu y^{d-\alpha}}+\eta_y(t),
\end{align}
where $\lim_{t\to\infty}\sup_{y\geq \rho}\eta_y(t)= 0$ for any $\rho>0$. The proof of~\cite[Prop.~3.8]{HMS06}
also contains the  idea about the right scaling for $Z_t$. Here we identify the limiting law, stated
in Theorem~\ref{thm-Ztlim}.

\begin{lemma} 
\label{psi}
As $t\to\infty$, the variable ${Z_t}/r_t$ converges weakly towards a random variable $X$ 
with Lebesgue density $p$, which was defined in Theorem~\ref{thm-Ztlim}.
\end{lemma}

\begin{Proof}{Proof.}
Let $A\subset \R^d$ be measurable with Leb$(\partial A)=0$. It suffices
to show that $$\lim_{t\to\infty}\Prob(Z_t/r_t\in A)=\int_{A} p(x) \, \dd x\, .$$
Let $\e>0$ and recall that $d-\alpha<0$. Pick $\rho>0$ so small that $\exp\{-\mu {\rho}^{d-\alpha}\}<\e/4$
and
\begin{align*}
\int_0^{\rho}\int_{A}\frac{\alpha\exp\left\{-\mu y^{d-\alpha}\right\}
}{(y+q|x|)^{\alpha+1}}\dd x\dd y<\e/4,
\end{align*}
which is possible since for $\rho=\infty$ the left hand side is equal to $\int_{A} p(x)\, \dd x \le 1$,
by Remark~\ref{3}.  Let $\eta_y(t)$ be as in \eqref{psiz}. 
Further, choose $T$ such that for all $t>T$ one has $\eta_y(t)<\e/4$ for all $y\in [\rho,\infty)$
and, moreover, 
\begin{align*}
\sup_{y\ge \rho}|\eta_y(t)|\int_{0}^{\infty}\int_{A}\frac{\alpha
\dd x\dd y}{(y+q|x|)^{\alpha+1}}<\e/4,
\end{align*}
which is possible as the integral is finite and $\eta_y(t)\to 0$ uniformly for $y\ge \rho$. 
We have
\begin{equation}
\label{decomp}
\Prob\Big(\frac{Z_t}{r_t}\in A\Big)
\leq \Prob\left(\Psi_t(Z_t)<{a}_t\rho\right)+
\Prob\left(Z_t\in r_t A,\Psi_t(Z_t)\ge{a}_t\rho\right).
\end{equation} 
For any $t>T$, the first probability can be bounded, with the help of \eqref{psiz}, by
\begin{align}
\label{eps2}
\Prob\left(\Psi_t(Z_t)<{a}_t\rho\right)
=\exp\big\{-\mu \rho^{d-\alpha}\big\}+\eta_{\rho}(t)
< \e/2.
\end{align}
Further, for any $t>T$, we compute the second probability as follows:
%\begin{align*}
%\Prob&\left(\frac{\Psi_t(Z_t)}{{a}_t}\le y,\frac{|Z_t|}{r_t}\in (a,b), 
%\frac{Z_t}{|Z_t|}\in A\right)
%=\sum_{z\in}\Prob\left(\Psi_t(Z_t)\le {a}_ty, Z_t=z\right)\\
%&=\sum_{z\in}\Prob\left(\Psi_t(z)\le {a}_ty, Z_t=z\right)
%=\sum_{z\in}\int_0^y\Prob\left(\Psi_t(z)\le {a}_tx\forall \widetilde{z}\neq z, 
%\Psi_t(z)={a}_tdx
%\right)
%\end{align*}
\begin{equation}
\label{1ugly}
\begin{aligned}
\Prob&\left(Z_t\in r_t A,\Psi_t(Z_t)\ge{a}_t\rho\right)
=\sum_{z\in r_t A\cap\Z^d}\int_{\rho}^{\infty}\Prob\left(Z_t=z,{a}_t^{-1}\Psi_t(Z_t)\in \dd y\right)
\notag\\
&
=\sum_{z\in r_t A\cap\Z^d}\int_{\rho}^{\infty}\Prob\left(\Psi_t(\widetilde{z})<{a}_ty~\forall\widetilde{z}\neq z,
{a}_t^{-1}\Psi_t(z)\in \dd y\right)\notag\\
&=\sum_{z\in r_t A\cap\Z^d}\int_{\rho}^{\infty}\Big[\prod_{\widetilde{z}\in \Z^d\setminus\{z\}}\Prob\Big(\xi(\widetilde{z})<{a}_ty+
\frac{|\widetilde{z}|}{t}\log\frac{|\widetilde{z}|}{2det}\Big)\Big]\\
&\qquad\qquad\times \Prob\Big(
{a}_t^{-1}\left[\xi(z)-\frac{|z|}{t}\log\frac{|z|}{2det}\right]\in \dd y\Big)\notag\\
%&=\sum_{z\in G_{t,A}^{b,c}}\int_{\rho}^{\infty}{a}_tF'\left(
%{a}_ty+\frac{|z|}{t}\log\frac{|z|}{2det}\right)\prod_{\widetilde{z}\neq z}
%F\left({a}_ty+
%\frac{|\widetilde{z}|}{t}\log\frac{|\widetilde{z}|}{2det}\right)dy\\
&=\int_{\rho}^{\infty}\!\!\!\!\dd y
\Big[\prod_{\widetilde{z}\in \Z^d}F\left({a}_ty+
\frac{|\widetilde{z}|}{t}\log\frac{|\widetilde{z}|}{2det}\right)\Big]
\sum_{z\in r_t A\cap\Z^d}
{a}_t
\frac{F'\left({a}_ty+\frac{|z|}{t}\log\frac{|z|}{2det}\right)}
{F\left({a}_ty+\frac{|z|}{t}\log\frac{|z|}{2det}\right)},
\end{aligned}
\end{equation}
since $\frac 1a (\xi(0)-K)$ has the density $y\mapsto a F'(ya+K)$. 
\smallskip

Recall \eqref{psiz} and $F(x)=1-x^{-\alpha}$ and therefore $a_t\frac{F'({a}_tv)}{F({a}_tv)}=
\frac{\alpha}{v({a}_t^{\alpha}v^{\alpha}-1)}$. Hence, with 
$$v_t(y,z)=y+\frac{|z|}{a_t t}\log\frac{|z|}{2det},$$
we obtain
$$
\begin{aligned}
\Prob&\left(Z_t\in r_t A,\Psi_t(Z_t)\ge{a}_t\rho\right)\\
&=\int_{\rho}^{\infty}\!\!\!\!\dd y\,\Big[{\rm e}^{-\mu y^{d-\alpha}}+\eta_y(t)\Big]\sum_{z\in r_t A\cap\Z^d}
\frac{\alpha}{v_t(y,z)({a}_t^{\alpha}v_t(y,z)^{\alpha}-1)}.
\end{aligned}
$$
Using that $a\log\frac a{2de}\geq -2d$ for any $a\geq0$, we have $v_t(y,z)\geq y-\frac {2d}{a_t}=y(1+o(1))$, uniformly for $y\geq \rho$. Hence,
\begin{equation}\label{Ztproof1}
\begin{aligned}
\Prob&\left(Z_t\in r_t A,\Psi_t(Z_t)\ge{a}_t\rho\right)\\
&=(1+o(1))\int_{\rho}^{\infty}\!\!\!\!\dd y\,\Big[{\rm e}^{-\mu y^{d-\alpha}}+\eta_y(t)\Big]
\sum_{z\in r_t A\cap\Z^d}\frac{\alpha}{{a}_t^{\alpha}v_t(y,z)^{\alpha+1}}.
\end{aligned}
\end{equation}
Fix some small $\delta>0$ and put $f_t=(\log t)^{-\delta}$ and $g_t=(\log t)^\delta$. We divide the sum over  $z\in r_t A\cap\Z^d$ on the 
right hand side of \eqref{1ugly} into the three parts where $|z|< r_tf_t$, $ r_tf_t\leq |z|\leq r_t g_t$ and $r_t g_t<|z|$. Hence, using 
an obvious notation in \eqref{Ztproof1},
\begin{equation}
\Prob\left(Z_t\in r_t A,\Psi_t(Z_t)\ge{a}_t\rho\right)
=I_t+II_t+III_t.
\end{equation}
We show that $I_t$ and $III_t$ vanish and that $|II_t-\int_A p(x)\dd x|\leq \e/2+o(1)$ as $t\to\infty$. Combining this with \eqref{decomp} 
and \eqref{eps2}, the convergence of $Z_t/r_t$ to the distribution with density $p$ follows.
\smallskip

We start with the estimate for $I_t$. From \eqref{Ztproof1} and $v_t(y,z)\geq y(1+o(1))$, we have
$$
I_t\leq \int_{\rho}^{\infty}\!\!\!\!\dd y\,\Big[{\rm e}^{-\mu y^{d-\alpha}}+\eta_y(t)\Big]\sum_{z\in\Z^d,|z|\leq r_tf_t}\frac{O(1)}
{a_t^\alpha y^{\alpha+1}}
\leq \int_{\rho}^{\infty}\!\!\!\!\dd y\,\Big[{\rm e}^{-\mu y^{d-\alpha}}+\eta_y(t)\Big]\frac{ O(1) f_t^d}{y^{\alpha+1}},
$$
where $O(1)$ does not depend on  $y$ nor on $z$, and we have used that $r_t^d=a_t^\alpha$. Since $\lim_{t\to\infty}f_t^d=0$, we see 
that $\lim_{t\to\infty} I_t=0$.

Now we turn to $II_t$. Recall that $q=\frac{d}{\alpha-d}$.
For $r_tf_t\leq |z|\leq r_t g_t$ we have 
\begin{equation}
\label{justlog}
\log\frac{|z|}{2det}=q\,(1+o(1))\,\log t.
\end{equation} 
Using \eqref{justlog} and the relations $t{a}_t=r_t\log t$
and $r_t^{d}={a}_t^{\alpha}$, we obtain, uniformly for $z\in r_t A\cap \Z^d$ satisfying $r_t f_t\leq|z|\leq r_t g_t$, and uniformly 
\mbox{for~$y\in (\rho,\infty)$},
$$
\begin{aligned}
\frac{\alpha}{{a}_t^{\alpha}v_t(y,z)^{\alpha+1}}
&=(\alpha+o(1))a_t^{-\alpha}\Big(y+\frac{|z|}{r_t\log t} q\,(1+o(1))\,\log t\Big)^{-\alpha-1}\\
&=(\alpha+o(1))r_t^{-d}\Big(y+\frac {|z|}{r_t} q\Big)^{-\alpha-1}.
\end{aligned}
$$
Substituting this into~\eqref{1ugly}, using~\eqref{psiz}, Leb$(\partial A)=0$, 
and interchanging the integrals gives
\begin{equation}
\label{last}
\begin{aligned}
II_t
&=(1+o(1))\, \int_{\rho}^{\infty}\left(e^{-\mu y^{d-\alpha}}+\eta_y(t)\right)
\int_{A}\1\{f_t\leq |x|\leq g_t\}\frac{\alpha}{(y+q|x|)^{\alpha+1}}\, \dd x \, \dd y\\
&=\int_{A}p(x)\dd x-
\int_0^{\rho}\int_{A}\frac{\alpha e^{-\mu y^{d-\alpha}}
\dd x\dd y}{(y+q|x|)^{\alpha+1}}
+\int_{\rho}^{\infty}\int_{A}\frac{\alpha\eta_y(t)
\dd x\dd y}{(y+q|x|)^{\alpha+1}}+o(1).
\end{aligned}
\end{equation}
Hence, by our choice of $\rho$, we have that $|II_t-\int_A p(x)\dd x|\leq \e/2+o(1)$.
\smallskip

Finally, we estimate $III_t$. For $|z|\geq r_tg_t$, we estimate $\log\frac{|z|}{2det}\geq \log \frac{r_tg_t}{2det}=(q+o(1))\log t$ 
and use the monotonicity to estimate, in the same way as for the term $II_t$,
$$
\begin{aligned}
III_t&\leq  O(1)\, \int_{\rho}^{\infty}\left(e^{-\mu y^{d-\alpha}}+\eta_y(t)\right)
\int_{\R^d}\1\{|x|\geq g_t\}\, \frac{1}{(y+q|x|)^{\alpha+1}}\, \dd x \, \dd y\\
&\leq O(1)\, \int_{\R^d}\1\{|x|\geq g_t\}\, p(x)\dd x.
\end{aligned}
$$
Since $p$ is integrable over $\R^d$ and $\lim_{t\to\infty}g_t=\infty$, we also have that $III_t$ vanishes as $t\to\infty$. This finishes the proof.
\qed\end{Proof}
\medskip

We now quantify the difference between the largest and the second-largest value of $\Psi_t$ in terms 
of their joint limit law. Recall the definition of $Z_t^{\ssup 1}$ and $Z_t^{\ssup 2}$ from the
beginning of this section, and also that $Z_t$ is identical to $Z_t^{\ssup 1}$.
\medskip

\begin{lemma} 
\label{psi2}
${{a}_t}^{-1}(\Psi_t(Z_t^{\ssup 1}),\Psi_t(Z_t^{\ssup 2}))\Longrightarrow
(Y_1,Y_2)$ weakly as $t\to\infty$, where $(Y_1,Y_2)$ is a $(0,\infty)\times (0,\infty)$-valued 
random variable with  distribution function
\begin{align*}
P(Y_1\le y_1,Y_2\le y_2)=\left\{
\begin{array}{ll}
\exp\left\{-\mu y_2^{d-\alpha}\right\}\left[1+\mu\left(y_2^{d-\alpha}-y_1^{d-\alpha}\right)\right]& 
\text{if }0< y_2\le y_1, \\[2mm]
\exp\left\{-\mu y_1^{d-\alpha}\right\}& \text{if }0 <  y_1< y_2.
\end{array}
\right.
\end{align*} 
\end{lemma}

\begin{Proof}{Proof.} First we argue that $\Psi_t(Z_t^{\ssup 1})\geq\Psi_t(Z_t^{\ssup 2})\geq 0$ almost surely,  
for all sufficiently large~$t\geq 0$. Indeed, since $\Psi_t(Z_t^{\ssup 2})$ is the second-largest of the values 
$\xi(z)-\frac{|z|}{t}\log\frac{|z|}{2det}$ with $z\in\Z^d$, it may be bounded from below against the minimum of any two of these. 
Picking $z=0$ and $z=z_0$ equal to a neighbour of the origin, we obtain the lower bound
\begin{align*}
\Psi_t(Z_t^{\ssup 2})&\ge \min\left\{\xi(0),\xi(z_0)-\sfrac{1}{t}\log\sfrac{1}{2det}\right\},
\end{align*}
which is nonnegative for all sufficiently large $t$. Hence it is sufficient to consider $y_1,y_2>0$. 
\smallskip

First, consider the case $0< y_1<y_2$. Using that $\Psi_t(Z_t^{\ssup 2})\le \Psi_t(Z_t^{\ssup 1})$
and~\eqref{psiz}, we obtain
\begin{align*}
\Prob\left(\Psi_t(Z_t^{\ssup 1})\le {a}_ty_1,\Psi_t(Z_t^{\ssup 2})\le {a}_ty_2\right)
=\Prob\left(\Psi_t(Z_t^{\ssup 1})\le {a}_ty_1\right)\to \exp\big\{-\mu y_1^{d-\alpha}\big\}.
\end{align*}

Second, assume $0< y_2\le y_1$. Using that $\Psi_t(Z_t^{\ssup 2})\le \Psi_t(Z_t^{\ssup 1})$, we obtain
\begin{equation}\label{prejoint}
\begin{aligned}
\Prob&\left(\Psi_t(Z_t^{\ssup 1})\le {a}_ty_1,\Psi_t(Z_t^{\ssup 2})\le {a}_ty_2\right)\\
&=\Prob\left(\Psi_t(Z_t^{\ssup 1})\le {a}_ty_2\right)+\Prob\left({a}_ty_2< \Psi_t(Z_t^{\ssup 1})\le {a}_ty_1,\Psi_t(Z_t^{\ssup 2})\le {a}_ty_2\right).
\end{aligned}
\end{equation}
Because of \eqref{psiz}, it suffices to study the second term on the right. Taking into account the independence of 
the random variables $(\Psi_t(z) \colon z\in\Z^d)$, we compute
we obtain 
\begin{equation}
\label{joint}
\begin{aligned}
\Prob&\left({a}_ty_2< \Psi_t(Z_t^{\ssup 1})\le {a}_ty_1,\Psi_t(Z_t^{\ssup 2})\le {a}_ty_2\right)\\
&=\sum_{z\in\Z^d}\Prob\left({a}_ty_2<\Psi_t(z)\le {a}_ty_1, \Psi_t(\widetilde{z})\le {a}_ty_2~
\forall\widetilde{z}\neq z\right)\\
&=\sum_{z\in\Z^d}\frac{\Prob\big({a}_ty_2<\Psi_t(z)\le {a}_ty_1\big)}
{\Prob\big(\Psi_t(z)\le {a}_ty_2\big)}
\prod_{\widetilde{z}\in\Z^d}\Prob\big(\Psi_t(\widetilde{z})\le {a}_ty_2\big)\\
&=\Prob\big(\Psi_t(Z_t^{\ssup 1})\le {a}_ty_2\big)
\sum_{z\in\Z^d}\frac{\Prob\big({a}_ty_2<\xi(z)-\frac{|z|}{t}\log\frac{|z|}{2det}\le {a}_ty_1\big)}
{\Prob\big(\xi(z)\le {a}_ty_2+\frac{|z|}{t}\log\frac{|z|}{2det}\big)}\\
%\end{aligned}
%\end{equation}
%\begin{equation}
%\begin{aligned}
&=\Prob\left(\Psi_t(Z_t^{\ssup 1})\le {a}_ty_2\right)\sum_{z\in\Z^d}\frac{\overline{F}\big({a}_ty_2+\frac{|z|}{t}\log\frac{|z|}{2det}\big)
-\overline{F}\big({a}_ty_1+\frac{|z|}{t}\log\frac{|z|}{2det}\big)}
{F\big({a}_ty_2+\frac{|z|}{t}\log\frac{|z|}{2det}\big)},
\end{aligned}
\end{equation}
where $\overline{F}(x)=1-F(x)=x^{-\alpha}$ is the tail of the distribution. Note that all the denominators are positive for all 
sufficiently large $t$.\smallskip
 
Since $a\log\frac{a}{2de}\geq -2d$ for any $a\geq 0$, we have $F({a}_ty_2+\frac{|z|}{t}\log\frac{|z|}{2det})\ge 
F({a}_ty_2-2d)=1+o(1)$ uniformly in $z$. To calculate the 
numerator, fix some small $\delta>0$ and denote $f_t=(\log t)^{-\delta}$ and
$g_t=(\log t)^{\delta}$. We split the sum into the three parts, as to where $|z|/r_t$ 
is smaller than $f_t$, between $f_t$ and $g_t$ and larger than $g_t$. We show next
that the two boundary contributions vanish, while the middle one has a nontrivial limit.
\smallskip

First, consider the domain where $|z|<r_tf_t$. We have 
\begin{align*}
\frac{|z|}{t}\log\frac{|z|}{2det}\le 
\frac{r_tf_t}{t}\log\frac{r_tf_t}{2det}=\frac{qf_tr_t\log t(1+o(1))}{t}=qf_t{a}_t(1+o(1))=o({a_t}),
\end{align*}
which, together with $r_t^d={a}_t^{\alpha}$, implies, for any $y>0$,
\begin{align}
\sum_{|z|<r_tf_t} \!\!\!\overline{F}\Big({a}_ty+\frac{|z|}{t}\log\frac{|z|}{2det}\Big)=
\sum_{|z|<r_tf_t} \!\!\!
\overline{F}\left({a}_ty(1+o(1))\right)=O\big((r_tf_t)^d\big)({a}_ty)^{-\alpha}=o(1).
\label{sum1}
\end{align}
Second, consider the domain where $r_tf_t\le |z|\le r_tg_t$. In this case $\log\frac{|z|}{2det}=q\log t (1+o(1))$ uniformly in $z$. 
Hence, using $\overline F(x)=x^{-\alpha}$, $r_t\log t=
t{a}_t$ and ${a}_t^{\alpha}=r_t^d$, we obtain  
\begin{equation}\label{aaa}
\overline{F}\Big({a}_ty+\frac{|z|}{t}\log\frac{|z|}{2det}\Big)
=\overline{F}\Big({a}_t\big(y+q\frac{|z|}{r_t}(1+o(1))\big)\Big)
=(1+o(1))\,r_t^{-d} \, \Big(y+q\frac{|z|}{r_t}\Big)^{-\alpha}.
\end{equation}
Summing over $r_tf_t\le |z|\le r_tg_t$, and turning the sum into an integral, we obtain
\begin{equation}\label{sum2}
\begin{aligned}
\sum_{f_t r_t\le |z|\le g_tr_t}\!\!\overline{F}\Big({a}_ty+\frac{|z|}{t}\log\frac{|z|}{2det}\Big)
&=(1+o(1))\, \int_{\R^d}\1\{f_t\leq |x|\leq g_t\}\big(y+q|x|\big)^{-\alpha}\,\dd x\\
&=(1+o(1))\, \mu \, y^{d-\alpha},
\end{aligned}
\end{equation}
where we use \cite[Lemma~3.9]{HMS06} to evaluate the integral and recall the definition of $\mu$ from~\eqref{mudef}. 

Finally, consider $|z|>r_tg_t$. Since  $\log\frac{|z|}{2det}\ge q\log t(1+o(1))$ 
uniformly in $z$, we have \lq$\leq$\rq\ instead of the first equality in \eqref{aaa}.
By the same procedure as in the case $r_tf_t\le |z|\le r_tg_t$, 
\begin{equation}\label{sum3}
\sum_{|z|>r_tg_t}\overline{F}\Big({a}_ty+\frac{|z|}{t}\log\frac{|z|}{2det}\Big)
\le (1+o(1))\int_{\R^d}\1\{|x|\geq g_t\}\big(y+q|x|\big)^{-\alpha}\,\dd x=o(1).
\end{equation}

Using~\eqref{sum1},~\eqref{sum2}, and~\eqref{sum3}, we obtain, for any $y>0$,
\begin{align*}
\sum_{z}\overline{F}\Big({a}_ty+\frac{|z|}{t}\log\frac{|z|}{2det}\Big)=\mu y^{d-\alpha}+o(1).
\end{align*}
Using this and \eqref{psiz} in~\eqref{joint} and substituting this and again \eqref{psiz} in \eqref{prejoint}, we obtain 
\begin{align*}
\Prob\left(\Psi_t(Z_t^{\ssup 1})\le {a}_ty_1,\Psi_t(Z_t^{\ssup 2})\le {a}_ty_2\right)
=\exp\big\{-\mu y_2^{d-\alpha}\big\}\,\left[1+\mu\big(y_2^{d-\alpha}-y_1^{d-\alpha}\big)\right]+o(1),
\end{align*}
which completes the proof.
\qed\end{Proof}

%\begin{lemma} 
%\label{xi1xi2}
%There exists a function $b:(0,\infty)\to (0,\infty)$ such that $b_r\to\infty$ as 
%$r\to\infty$ and $\Prob\left(\xi_r^{\ssup 1}-\xi_r^{\ssup 2}\ge b_r\right)\to 1$. 
%\end{lemma}

%\begin{proof} 
%It follows from 

%We have
%\begin{align*}
%\Prob\left(\xi_r^{\ssup 1}-\xi_r^{\ssup 2}> b_r\right)
%&=\kappa_dr^d(1+o(1))\Prob\left(\xi(0)-
%\max_{0<|z|\le r}\xi(z)>b_r\right)\\
%&=\kappa_dr^d(1+o(1))\int_{1}^{\infty} F(x)^{\kappa_dr^d(1+o(1))}dF(x+b_r)
%\end{align*}
%\end{proof}

\section{Proof of Theorem~\ref{main}: Complete localisation}\label{sec-proofmain}

In this section, we prove Theorem~\ref{main}. Section~\ref{sec-uDecomp} presents the details of the
decomposition of $u$ into three parts, which is announced informally in Section~\ref{overview}.
Subject to the two main propositions, whose proofs are deferred to Section~\ref{sec-u1u2esti} and Section~\ref{sec-u3esti}, 
we finish the proof of Theorem~\ref{main} in this section. 
Proposition~\ref{lu1u2} is proved in Section~\ref{sec-u1u2esti}, where we show that the total mass of 
the first two contributions is negligible, using extreme value theory and certain limit laws. 
Proposition~\ref{lu3} is proved in Section~\ref{sec-u3esti}, where we show that the third contribution 
is asymptotically concentrated in $Z_t$.

\subsection{Decomposing $\boldsymbol u$.}\label{sec-uDecomp}\\[-2mm]

Let $(X_s \colon s\in [0,\infty))$ be the continuous-time simple random walk on $\Z^d$ with generator 
$\Deltad$. By $\P_z$ and $\E_z$ we denote the probability measure and the expectation with respect to 
the walk starting at $z\in\Z^d$. According to \cite[Theorem~2.1]{GM90}, the unique nonnegative
solution of~\eqref{pam} can  be expressed in terms of the Feynman-Kac formula as 
\begin{align}
\label{FC1}
u(t,z)=\E_0 \Big[\exp\Big\{\int_0^t\xi(X_s)\dd s\Big\}\one\{X_t=z\}\Big], \qquad t>0, \, z\in\Z^d,
\end{align} 
where we also used the time-reversal property of the random walk. 
We denote the entrance time into a set $A\subset \Z^d$ by $\tau_A=\inf\left\{t\ge 0\colon X_t\in A\right\}$ 
and abbreviate $\tau_z=\tau_{\left\{z\right\}}$.
By $B_R=\{z\in \Z^d\colon |z|\leq R\}$ we denote the box in $\Z^d$ with radius $R>0$.
Let $h\colon (0,\infty)\to (0,\infty)$ be such that 
\begin{equation}
\label{h}
\lim_{t\to\infty}h_t=0\qquad\text{and}\qquad\lim_{t\to\infty}h_t\frac{\sqrt{\log t} }{\log\log t}=\infty,
\end{equation}
and define the \emph{random} radius 
\begin{equation}\label{Rtdef}
R_t=|Z_t|(1+h_t).\\[2.5mm]
\end{equation}
We write $u(\theta,z)=u_1^{\ssup t}(\theta,z)+u_2^{\ssup t}(\theta,z)+u_3^{\ssup t}(\theta,z)$,
where 
\begin{align*}
u_1^{\ssup t}(\theta,z)&=\E_0\Big[\exp\Big\{\int_0^{\theta}\xi(X_{s})\dd s\Big\}\1\{X_{\theta}=z\}\,
\1\{\tau_{B_{R_t}^{\mathrm{c}}}\le \theta\}\Big]\\
u_2^{\ssup t}(\theta,z)&=\E_0\Big[\exp\Big\{\int_0^{\theta}\xi(X_{s})\dd s\Big\}\1\{X_{\theta}=z\}\,
\1\{\tau_{B_{R_t}^{\mathrm{c}}}>\theta\}\1\{\tau_{Z_t}>\theta\}\Big]\\
u_3^{\ssup t}(\theta,z)&=\E_0\Big[\exp\Big\{\int_0^{\theta}\xi(X_{s})\dd s\Big\}\1\{X_{\theta}=z\}\,
\1\{\tau_{B_{R_t}^{\mathrm{c}}}>\theta\}\1\{\tau_{Z_t}\le \theta\}\Big],
\end{align*}
for $(\theta,z)\in (0,\infty)\times \Z^d$ and $t>0$. We are mainly interested in this decomposition for 
$\theta=t$.\medskip

\begin{prop}[Estimating ${u_1^{\ssup t}}$ and ${u_2^{\ssup t}}$]\label{lu1u2}
$$
\lim_{t\to\infty}\frac{\sum_{z\in\Z^d} u_1^{\ssup t}(t,z)}{U(t)}= 0\qquad\mbox{and}\qquad
\lim_{t\to\infty}\frac{\sum_{z\in\Z^d} u_2^{\ssup t}(t,z)}{U(t)}= 0\qquad \mbox{in probability.}
$$
\end{prop}\medskip

\begin{prop}[Estimating ${u_3^{\ssup t}}$]\label{lu3}
$$
\lim_{t\to\infty}\frac{\sum_{z\in\Z^d\setminus\{Z_t\}} u_3^{\ssup t}(t,z)}{U(t)}= 0\qquad\mbox{in probability.}
$$
\end{prop}\medskip

These two propositions will be proved in the next two sections.
Using them, we can easily finish the proof of our first main result:
\medskip

\begin{Proof}{Proof of Theorem~\ref{main}.} Recall that $u=u_1^{\ssup t}+u_2^{\ssup t}+u_3^{\ssup t}$. Since
$u_1^{\ssup t}$ and $u_2^{\ssup t}$ are nonnegative, %we obtain
\begin{align*}
\frac{\sum_{z\in\Z^d\setminus\{Z_t\}} u(t,z)}{U(t)}
\leq \frac{\sum_{z\in\Z^d} u_1^{\ssup t}(t,z)}{U(t)}+\frac{\sum_{z\in\Z^d} u_2^{\ssup t}(t,z)}{U(t)}
+\frac{\sum_{z\in\Z^d\setminus\{Z_t\}} u_3^{\ssup t}(t,z)}{U(t)},
\end{align*}
and the right hand side vanishes in probability as $t\to\infty$, by Propositions~\ref{lu1u2} and~\ref{lu3}.
\qed\end{Proof}
\bigskip

\subsection{Proof of Proposition~\ref{lu1u2}:
Estimating $\boldsymbol{u_1^{\ssup t}}$ and $\boldsymbol{u_2^{\ssup t}}$.}\label{sec-u1u2esti}
\\[-2mm]

In this section we prove Proposition~\ref{lu1u2}, i.e., we show that the contributions coming from $u_1^{\ssup t}$ 
and $u_2^{\ssup t}$ are negligible. To prepare this, we first show that $Z_t=Z_t^{\ssup 1}$, the maximal point of~$\Psi_t$, 
is also maximal for the potential~$\xi$ in the smallest centred box that contains it. Then we show that, by our choice of $R_t$ 
in \eqref{Rtdef}, $Z_t$ is also maximal for $\xi$ in the box with radius~$R_t$. Finally, the difference to the second-largest 
value of $\xi$ in this box diverges. 
\smallskip

In order to formulate these statements, we define the two upper order statistics for the potential $\xi$ by 
\begin{align*}
\xi_r^{\ssup 1}=\max\left\{\xi(z)\colon |z|\le r\right\}\qquad\text{and}\qquad
\xi_r^{\ssup 2}=\max\left\{\xi(z)\colon |z|\le r, \xi(z)\neq \xi_r^{\ssup 1}\right\}.
\end{align*}
It follows from the continuity of distribution of $\xi(0)$ that, for any $r>0$, each of the sets
$\left\{x\in\Z^d\colon |x|\leq r, \xi(x)=\xi_r^{\ssup i}\right\}$, $i=1,2$, contains exactly one point, almost surely. 

\begin{lemma} \ \\[-5mm]
\label{xi}
\begin{enumerate}
\item[(i)] $\displaystyle\lim_{t\to\infty}\Prob\big(\xi(Z_t)=\xi_{|Z_t|}^{\ssup 1}\big)= 1$;\\[-2mm]

\item[(ii)] $\displaystyle\lim_{t\to\infty}\Prob\big(\xi(Z_t)=\xi_{R_t}^{\ssup 1}\big)= 1$;\\[-2mm]

\item[(iii)] $\displaystyle\lim_{t\to\infty}\Prob\big(t\xi(Z_t)>|Z_t|\big)= 1$.\\[-2mm]

\item[(iv)] There exists ${b}\colon (0,\infty)\to (0,\infty)$ such that $\lim_{t\to\infty}{b}_t=\infty$
and $$\lim_{t\to\infty}\Prob\big(\xi_{R_t}^{\ssup 1}- \xi_{R_t}^{\ssup 2}\ge {b}_t\big)= 1.$$ 
\end{enumerate}
\end{lemma}

\begin{Proof}{Proof.}
(i) Set $f_t=(\log t)^{-\delta}$ for some small $\delta>0$. By Lemma~\ref{psi},  
$\lim\Prob\left(|Z_t|\ge f_t r_t\right)=1$. It thus suffices to show
that $\xi(Z_t)=\xi_{|Z_t|}^{\ssup 1}$ on the set 
$\left\{|Z_t|\ge f_t r_t\right\}$ for all large $t$.

Suppose for contradiction that $|Z_t|\ge f_t r_t$, but there exists $z\neq Z_t$ such that $|z|\le |Z_t|$ 
and $\xi(z)>\xi(Z_t)$. Since $f_t r_t/t\to \infty$, we may assume
that $t$ is large enough to satisfy $\frac{f_t r_t}{t}\,\log(\frac{f_t r_t}{2det})>0$.   
Using that $r\mapsto \frac{r}{t}\log\frac{r}{2det}$ is increasing for $r> 2det$ and nonpositive
otherwise, we get 
\begin{align*}
\Psi_t(z)=\xi(z)-\frac{|z|}{t}\log\frac{|z|}{2det}>\xi(Z_t)-
\frac{|Z_t|}{t}\log\frac{|Z_t|}{2det}=\Psi_t(Z_t),
\end{align*}
which contradicts the fact that $Z_t$ is the maximum of $\Psi_t$. Hence (i) is proved.
\pagebreak[3]\smallskip

(ii) As $R_t\ge |Z_t|$ we clearly have $\xi_{R_t}^{\ssup 1}\ge \xi_{|Z_t|}^{\ssup 1}$. 
Let $f,g\colon (0,\infty)\to (0,\infty)$ be such that $f_t\to 0$, $g_t\to\infty$ and 
$h_t\log\left(g_t/f_t\right)\to 0$. Then we obtain
\begin{equation}\label{ximaxproof1}
\begin{aligned}
\Prob&\big(\xi_{R_t}^{\ssup 1}>\xi_{|Z_t|}^{\ssup 1},r_tf_t\le |Z_t|\le r_tg_t,
\xi(Z_t)
=\xi_{|Z_t|}^{\ssup 1}\big)\\
&=\sum_{r=\lceil r_tf_t\rceil}^{\lfloor r_tg_t \rfloor}\Prob\big(\xi_{R_t}^{\ssup 1}>\xi_{|Z_t|}^{\ssup 1},|Z_t|=r,
\xi(Z_t)=\xi_{|Z_t|}^{\ssup 1}\big)\\
&\le \sum_{r=\lceil r_tf_t\rceil}^{\lfloor r_tg_t \rfloor}\Prob\big(\xi_{r(1+h_t)}^{\ssup 1}>\xi_r^{\ssup 1}>\xi_{r-1}^{\ssup 1}\big)\\
&=\sum_{r=\lceil r_tf_t\rceil}^{\lfloor r_tg_t \rfloor}\big[\Prob\big(\xi_r^{\ssup 1}>\xi_{r-1}^{\ssup 1}\big)
-\Prob\big(\xi_{r(1+h_t)}^{\ssup 1}=\xi_r^{\ssup 1}>\xi_{r-1}^{\ssup 1}\big)\big].
\end{aligned}
\end{equation}
Observe that, for any two finite non-empty subsets $A\subset B$ of $\Z^d$, we have $\Prob(\max_{z\in B}\xi(z)=\max_{z\in A}\xi(z))=|A|/|B|$, 
since all the values $\xi(z)$ with $z\in B$ are different, and the index of the maximal value is uniformly distributed over $B$. 
Also, observe that 
$$\big\{\xi_r^{\ssup 1}>\xi_{r-1}^{\ssup 1}\big\}=\big\{\max_{z\in B_r}\xi(z)=\max_{z\in \partial B_r}\xi(z)\big\}$$ and  $$\big\{\xi_{r(1+h_t)}^{\ssup 1}=\xi_r^{\ssup 1}>\xi_{r-1}^{\ssup 1}\big\}=\big\{\max_{z\in B_{r(1+h_t)}}\xi(z)=\max_{z\in \partial B_r}\xi(z)\big\},$$ where 
the \emph{inner boundary} of~$B_r$ is defined by 
$$\partial B_r = \big\{ x\in B_r \colon \mbox{ there is } y\not\in B_r \mbox{ with } |y-x|=1 \big\}\, .$$ 
Denoting $\sigma_d=\lim_{r\to\infty}|\partial B_r|r^{1-d}>0$ and 
$\kappa_d=\lim_{r\to\infty}|B_r|r^{-d}>0$, we therefore obtain from \eqref{ximaxproof1} that
\begin{equation}
\begin{aligned}
\Prob\Big(\xi_{R_t}^{\ssup 1} & >\xi_{|Z_t|}^{\ssup 1},r_tf_t\le |Z_t|\le r_tg_t,
\xi(Z_t)
=\xi_{|Z_t|}^{\ssup 1}\Big)\\
&\leq (1+o(1))\sum_{r=\lceil r_tf_t\rceil}^{\lfloor r_tg_t \rfloor}\Big[\frac{\sigma_d r^{d-1}}{\kappa_d r^d}
-\frac{\sigma_d r^{d-1}}{\kappa_d r^d(1+h_t)^d}\Big]\\
&=(1+o(1))\,\frac{\sigma_d}{\kappa_d}\, \Big[1
-\frac{1}{(1+h_t)^d}\Big]\, \sum_{r=\lceil r_tf_t\rceil}^{\lfloor r_tg_t \rfloor}\frac{1}{r}
=O(1)\,d \, h_t\, \log\frac{g_t}{f_t},
\end{aligned}
\end{equation}
and this vanishes as $t\to\infty$ because of our assumption on $f_t $ and $g_t$. Since we know from  Lemma~\ref{psi}, 
respectively from  (i), that the probabilities of the events $\{r_tf_t\le |Z_t|\le r_tg_t\}$ and 
$\{\xi(Z_t)=\xi_{|Z_t|}^{\ssup 1}\}$ tend to one as $t\to\infty$, the assertion (ii) is proved.
\smallskip

(iii) Let $f_t=1/\log t$. Using that $t\xi(Z_t)=\Psi_t(Z_t)+|Z_t|\log\frac{|Z_t|}{2det}$ and
$\Psi_t(Z_t)>0$, we obtain
\begin{align*}
\Prob&\big(t\xi(Z_t)\le|Z_t|, |Z_t|\ge r_tf_t\big)
=\Prob\big(t\Psi_t(Z_t)+|Z_t|\log\sfrac{|Z_t|}{2det}\le|Z_t|, 
|Z_t|\ge r_tf_t\big)\\
&\phantom{aaaaaaaaaaa}
\le \Prob\big(\log\sfrac{|Z_t|}{2det}\le 1, |Z_t|\ge r_tf_t\big)
\le \Prob\big(\log\sfrac{r_tf_t}{2det}\le 1\big).
\end{align*}
The right hand side is equal to zero
if $t$ is sufficiently large, as $r_tf_t/2det\to\infty$. By Lemma~\ref{psi}, the probability of the event $\{|Z_t|\ge r_tf_t\}$ 
tends to one, and this ends the proof of (iii).
\medskip\pagebreak[2]

(iv) There exists a scale function $\overline b\colon (0,\infty)\to(0,\infty)$ such that $\overline b_r\to\infty$ and 
$\Prob(\xi_r^{\ssup 1}-\xi_r^{\ssup 2}\geq \overline b_r)\to 1$ as $r\to\infty$. Indeed, this follows from the fact that 
the top two values of the order statistics satisfy the limit law
\begin{equation}
r^{-\frac{d}{\alpha}}\left(\xi_r^{\ssup 1},\xi_r^{\ssup 2}\right)\Longrightarrow (\Xi_1,\Xi_2)\qquad\mbox{as }r\to\infty,
\end{equation}
where $\Xi_1$ and $\Xi_2$ are two continuous $(0,\infty)$-valued random variables that satisfy $\Xi_1>\Xi_2$ almost surely, 
see \cite[Th.~4.2.8]{yellow} for the general limit assertion and \cite[p.~153]{yellow} for the discussion of the Pareto case. 
From this limit law, it is easy to construct the desired scale function $\overline b$.
Note that~(iv), which we now prove, does not follow immediately from this, because the radius~$R_t$ is 
chosen randomly. Define 
\begin{align*}
p_r=\Prob\left(\xi_r^{\ssup 1}-\xi_r^{\ssup 2}<\overline b_r\right).
\end{align*}
and choose $\overline{f}\colon (0,\infty)\to (0,\infty)$ in such a way that 
$\overline{f}_t\to 0$ and $r_t\overline{f_t}\to\infty$. 
As $p_r\to 0$ this implies $\overline{p}_t=\sup_{r>r_t\overline{f_t}}p_r\to 0$. 
Now we can choose $f\colon (0,\infty)\to (0,\infty)$ so that 
$f_t\to 0$, $f_t>\overline{f}_t$ and $\overline{p}_t\log f_t\to 0$.
Finally, we choose $g\colon (0,\infty)\to (0,\infty)$ 
such that $g_t\to\infty$ and $\overline{p}_t\log g_t\to 0$. This gives
\begin{align}
\label{fgc} 
\sup_{r\ge r_tf_t}p_{r(1+h_t)}\log\frac{g_t}{f_t}
\le\sup_{r\ge r_t\overline{f}_t}p_{r(1+h_t)}\log\frac{g_t}{f_t}
\le \overline{p}_t\log\frac{g_t}{f_t}\to 0.
\end{align}

Define
\begin{align*}
{b}_t=\inf_{r_tf_t\le r\le r_tg_t}\overline b_{r(1+h_t)}
\end{align*}
and note that ${b}_t\to \infty$ since $\overline b_r\to\infty$ and $r_tf_t\to \infty$.
Using the spatial homogeneity of the family $(\xi(z) \colon z\in\Z^d)$ and~\eqref{fgc}, we obtain
\begin{equation}
\begin{aligned}
\Prob&\left(\xi_{R_t}^{\ssup 1}-\xi_{R_t}^{\ssup 2}<{b}_t,r_tf_t\le |Z_t|\le r_tg_t,
\xi(Z_t)=\xi_{|Z_t|}^{\ssup 1}=\xi_{R_t}^{\ssup 1}\right)\\
&=\sum_{r_tf_t \le |z| \le r_tg_t}
\Prob\left(\xi_{R_t}^{\ssup 1}-\xi_{R_t}^{\ssup 2}<{b}_t,Z_t=z,
\xi(Z_t)=\xi_{|Z_t|}^{\ssup 1}=\xi_{R_t}^{\ssup 1}\right)\\
&\le \sum_{r_tf_t \le |z| \le r_tg_t}
\Prob\left(\xi_{|z|(1+h_t)}^{\ssup 1}-\xi_{|z|(1+h_t)}^{\ssup 2}<{b}_t,\xi(z)=\xi_{|z|(1+h_t)}^{\ssup 1}\right).
\end{aligned}
\end{equation}
Observe that the top two values in the order statistics are independent of the indices at which they are realised, i.e., the 
two events on the right hand side are independent, and that the probability of the second event is $1/|B_{|z|(1+h_t)}|$. 
As before, we denote $\sigma_d=\lim_{r\to\infty}|\partial B_r|r^{1-d}$ and $\kappa_d=\lim_{r\to\infty}|B_r|r^{-d}$. 
Using ${b}_t\leq \overline b_{|z|(1+h_t)}$ and the definition of $p_{r(1+h_t)}$, we therefore obtain 
\begin{equation}
\begin{aligned}
\Prob&\Big(\xi_{R_t}^{\ssup 1}-\xi_{R_t}^{\ssup 2}<{b}_t,r_tf_t\le |Z_t|\le r_tg_t,
\xi(Z_t)=\xi_{|Z_t|}^{\ssup 1}=\xi_{R_t}^{\ssup 1}\Big)\\
&\leq \sum_{r_tf_t \le |z| \le r_tg_t}
\Prob\Big(\xi_{|z|(1+h_t)}^{\ssup 1}-\xi_{|z|(1+h_t)}^{\ssup 2}<{b}_t\Big)
\,\Prob\Big(\xi(z)=\xi_{|z|(1+h_t)}^{\ssup 1}\Big)\\
&\le(1+o(1))\sum_{r_tf_t \le |z| \le r_tg_t} 
\Prob\Big(\xi_{|z|(1+h_t)}^{\ssup 1}-\xi_{|z|(1+h_t)}^{\ssup 2}<\overline b_{|z|(1+h_t)}\Big) \, \frac{1}{\kappa_d|z|^{d}(1+h_t)^d}\\
&\le O(1)\sum_{r=\lceil r_tf_t\rceil}^{\lfloor r_tg_t \rfloor} \frac{1}{r}\,
\sup_{r_tf_t\le r\le r_tg_t}p_{r(1+h_t)}
\le O(1)\,\log\frac{g_t}{f_t}\,\,\sup_{r\ge r_tf_t}p_{r(1+h_t)},\end{aligned}
\end{equation}
where we changed the sum over $z$ into a sum over $r$, which turns the term $|z|^{-d}$ into $\frac 1r$. By our assumptions, 
the right hand side vanishes as $t\to\infty$. Since we know from Lemma~\ref{psi} and~(i) and~(ii) that the probabilities of 
the events $\{r_tf_t\le |Z_t|\le r_tg_t\}$ and $\{\xi(Z_t)=\xi_{|Z_t|}^{\ssup 1}=\xi_{R_t}^{\ssup 1}\}$ tend to one, the proof 
of~(iv) is finished.
\qed\end{Proof}
\smallskip

Now we give a lower bound for the total mass $U(t)$ in terms of the maximal point $Z_t$ of $\Psi_t$. Recall that $O(t)$ 
denotes some deterministic function $(0,\infty)\to(0,\infty)$ that is at most linear in $t$ at infinity.

\begin{lemma}[Lower bound for $U(t)$]
\label{bounds}
\begin{align*}
\lim_{t\to\infty}\Prob\big(\log U(t)\ge t\xi(Z_t)-|Z_t|\log\xi(Z_t)+O(t)\big)= 1.
\end{align*}
\end{lemma}

\begin{Proof}{Proof.} 
Fix $\e\in(0,1-\frac d\alpha)$ and note from \eqref{asas} that $\xi_r^{\ssup 1}$ is asymptotically sublinear in $r$, almost surely. 
Hence, by~\cite[Lemma~2.2]{HMS06} there exists a random time $T$ such  that, for all $t>T$,
\begin{align}
\label{lb1}
\log U(t)\ge t\max_{0<\rho<1}\max_{z\in\Z^d}\left[(1-\rho)\xi(z)-\frac{|z|}{t}\log\frac{|z|}{e\rho t}\right]+O(t).
\end{align}  

On the event $\{t\xi(Z_t)>|Z_t|\}$, we substitute $\rho=\frac{|Z_t|}{t\xi(Z_t)}\in(0,1)$ and $z=Z_t$ in~\eqref{lb1} and obtain
\begin{align*}
\log U(t)\ge \big(t-\sfrac{|Z_t|}{\xi(Z_t)}\big)\xi(Z_t)-|Z_t|\log\sfrac{|Z_t|\xi(Z_t)}{e|Z_t|}+O(t)
=t \xi(Z_t)-|Z_t|\log\xi(Z_t)+O(t).
\end{align*}
By Lemma~\ref{xi}(iii) the probability of $\{t\xi(Z_t)>|Z_t|\}$ tends to one,
which implies the claim.
\qed\end{Proof}
\smallskip

Now we derive upper bounds for the total mass in terms of the sites 
$Z_t=Z_t^{\ssup 1}$ and $Z_t^{\ssup 2}$.

\begin{lemma}[Upper bounds for $u_1^{\ssup t}$ and $u_2^{\ssup t}$]\label{bounds2} \ \\ \vspace{-0.5cm}
\begin{enumerate}
\item[(i)]
$$\lim_{t\to\infty}\Prob\Big(\log \sum_{z\in\Z^d} u_2^{\ssup t}(t,z)\le 
t\Psi_t(Z_t^{\ssup 2})+O(t)\Big)= 1.$$
\item[(ii)]
$$\lim_{t\to\infty}\Prob\Big(\log \sum_{z\in\Z^d} u_1^{\ssup t}(t,z)
\le\max\Big\{ t\Psi_t(Z_t^{\ssup 2}),
t\xi(Z_t^{\ssup 1})-R_t\log\frac{R_t}{2det}\Big\}
+O(t)\Big)= 1.$$
\end{enumerate}
\end{lemma}

\begin{Proof}{Proof.}
(i) Note that 
$$
\sum_{z\in\Z^d} u_2^{\ssup t}(t,z)\leq 
\E_0\left[\exp\left\{\int_0^t\xi(X_s)\dd s\right\}\1\{\tau_{Z_t^{\ssup 1}}>t\}\right].
$$
Denote
\begin{align*}
\zeta_r=\max_{z\in B_r\setminus{ Z_t^{\ssup 1}}}\xi(z). 
\end{align*}
Denote by $J_t$ the number of jumps of the random walk $(X_s \colon s\ge 0)$ before time $t$. Note 
that $J_t$ has a Poisson distribution with parameter $2dt$, and that the path stays inside the box $B_{J_t}$ 
up to time~$t$. Therefore, on the event $\{\tau_{Z_t^{\ssup 1}}>t\}$, we can estimate $\xi(X_s)\le \zeta_{J_t}$ for $s\in[0,t]$. 
Summing over all values of $J_t$, we obtain
\begin{equation}\label{u2esti1}
\sum_{z\in\Z^d} u_2^{\ssup t}(t,z)
\leq \sum_{r=0}^{\infty}\E_0\left[\exp\left\{t\zeta_{J_t}\right\}\1\{\tau_{Z_t^{\ssup 1}}>t\}
\1\{J_t=r\}\right]\le \sum_{r=0}^{\infty}e^{t\zeta_r-2dt}\frac{(2dt)^r}{r!}.
\end{equation}
We now give an upper bound for the tail of the series on the right. Fix some $\theta>1$, $\eps>0$,
$1>\eta>d/\alpha$ and let $\beta=(1-\eta)^{-1}(1+\eps)$. Using Stirling's formula,
$$r!=\sqrt{2\pi r}\left(\sfrac{r}{e}\right)^re^{\delta(r)}, \qquad\mbox{ with }\lim_{r\uparrow\infty}\delta(r)=0\, ,$$
and the bound $\zeta_r \le \xi_r^{\ssup 1}\le r^{\eta}$ for all large~$r$,
we obtain, for all $r>t^{\beta}$ and large $t$, that
\begin{align*}
t&\,\zeta_r-2dt+r\log (2dt)-\log (r!)
\le tr^{\eta}-r\log\sfrac{r}{2d\ee t}-\delta(r)\\
&\le tr^{\eta}\left(1-\sfrac{r^{1-\eta}}{t}\log\sfrac{r}{2d\ee t}-\sfrac{\delta(r)}{tr^{\eta}}\right)
\le tr^{\eta}\left(1-t^{\eps}\log\sfrac{t^{\beta-1}}{2de}-\sfrac{\delta(r)}{tr^{\eta}}\right)\le -\theta\,\log r\, .
\end{align*}
Splitting the sum on the right of \eqref{u2esti1} at $r=\lceil t^\beta \rceil$
and noting that $\sum_{r>\lceil t^\beta \rceil} r^{-\theta}=o(1)$, we obtain
\begin{equation}\begin{aligned}
\label{est}
\log \sum_{z\in\Z^d} u_2^{\ssup t}(t,z) 
& \le  t\, \max_{0\le r\le t^{\beta}}\left[\zeta_r-\sfrac{r}{t}\log\sfrac{r}{2det}-\sfrac{1}{t}\log\sqrt{2\pi r}-
\sfrac{\delta(r)}{t}\right] -2dt +o(t) \\
& \le t\max_{r\in\N_0}\left[\zeta_r-\sfrac{r}{t}\log\sfrac{r}{2det}\right]+O(t).
\end{aligned}\end{equation}
Our goal is to show that the maximum on the right hand side is not larger than $\Psi_t(Z_t^{\ssup 2})$, the second-largest 
value of $\Psi_t$, with probability tending to one. \smallskip

Denote by $\rho_t$ the value at which this maximum is attained and let $z_t\in B_{\rho_t}\setminus\{Z_t^{\ssup 1}\}$ 
be the maximal point in the definition of $\zeta_{\rho_t}$. Both $\rho_t$ and $z_t$ are unique by the continuity of the potential distribution. 
For any $z\in\Z^d\setminus\{Z_t^{\ssup 1}\}$, we have, using the definition of $\Psi_t$, then the definition of 
$\zeta_{ |z|}$ and $z\neq Z_t^{\ssup 1}$, the definition of $\rho_t$ and finally the definition of~$z_t$,
\begin{equation}\label{Psitesti}
\begin{aligned}
\Psi_t(z)&=\xi(z)-\frac{|z|}{t}\log\frac{|z|}{2det}\le\zeta_{|z|}-\frac{|z|}{t}\log\frac{|z|}{2det}
\le \zeta_{\rho_t}-\frac{\rho_t}{t}\log\frac{\rho_t}{2det}\\
&=\xi(z_t)-\frac{\rho_t}{t}\log\frac{\rho_t}{2det}.
\end{aligned}
\end{equation}
On the event $\{\rho_t\geq 2det\}$, one can estimate $-\frac {\rho_t}t\log\frac{\rho_t}{2det}\leq -\frac {|z_t|}t\log\frac{|z_t|}{2det}$, 
since $|z_t|\leq\rho_t$, and since the map $r\mapsto \frac rt\log\frac{r}{2det}$ is increasing on $[2dt,\infty)$ and positive precisely 
on $(2det,\infty)$. Then \eqref{Psitesti} implies that $\Psi_t(z)\leq \Psi_t(z_t)$. Hence, $\Psi_t(z_t)$ turns out to be the second-largest 
value of $\Psi_t$, and it follows that $z_t=Z_t^{\ssup 2}$. Note that the last two terms of \eqref{Psitesti} are equal to the maximum on 
the right hand side of \eqref{est}, which therefore is not smaller than $\Psi_t(Z_t^{\ssup 2})$. Summarising, on the event $\{\rho_t\geq 2det\}$, 
we have the desired estimate. Hence, it remains to show that the probability of this event tends to one.
\smallskip

Recall that $q=d/(\alpha-d)$ and pick $\eps_1\in(0,q)$. It suffices to show that the probability of the event
$\{\rho_t<t^{q+1-\e_1}\}$ vanishes. For this purpose, pick $0<\e_2<\e_3<d/\alpha\e_1$ and $\e_4>0$ such that 
$\e_4(q+1-\e_1)<{d}/{\alpha}\e_1-\e_3$, and observe that
\begin{align}
\Prob\left(\rho_t<t^{q+1-\e_1}\right)
&\le \Prob\left(\Psi_t(Z_t^{\ssup 2})< t^{q-\e_2}\right)
+\Prob\left(\xi_{t^{q+1-\e_1}}^{\ssup 1}> t^{q-\e_3}\right)\notag\\
&\qquad +\Prob\left(\rho_t<t^{q+1-\e_1},
\Psi_t(Z_t^{\ssup 2})\ge t^{q-\e_2},\xi_{t^{q+1-\e_1}}^{\ssup 1}\le t^{q-\e_3}\right). 
\label{rhos}
\end{align}
The first term on the right hand side vanishes since, by Lemma~\ref{psi2}, $\Psi_t(Z_t^{\ssup 2})$ is of order 
$a_t=(t/\log t)^q$. The second term vanishes by \eqref{asas} applied to $t^{q+1-\e_1}$, because, 
almost surely, for any sufficiently large $t$,
\begin{align*}
\xi_{t^{q+1-\e_1}}^{\ssup 1}\le \left(t^{q+1-\e_1}\right)^{\frac{d}{\alpha}+\e_4} \le t^{q-\e_3}.
\end{align*}
Finally, we show that the third term is equal to zero for any sufficiently large $t$. Indeed, first estimate
\begin{equation}
\begin{aligned}
\Psi_t(Z_t^{\ssup 2})&=\xi(Z_t^{\ssup 2})-\sfrac{|Z_t^{\ssup 2}|}{t}\log\sfrac{|Z_t^{\ssup 2}|}{2det}
\leq\zeta_{|Z_t^{\ssup 2}|}-\sfrac{|Z_t^{\ssup 2}|}{t}\log\sfrac{|Z_t^{\ssup 2}|}{2det}\\
& \leq\max_{r\in\N_0}\left[\zeta_r-\sfrac{r}{t}\log\sfrac{r}{2det}\right]
=\zeta_{\rho_t}-\sfrac{\rho_t}{t}\log\sfrac{\rho_t}{2det}
\leq \xi_{\rho_t}^{\ssup 1}+2d,
\end{aligned}
\end{equation}
since $\frac rt\log\frac{r}{2det}\geq -2d$ for any $r\geq 0$. Hence, on the event 
$$\big\{\rho_t<t^{q+1-\e_1}, \, 
\Psi_t(Z_t^{\ssup 2})\ge t^{q-\e_2},\, \xi_{t^{q+1-\e_1}}^{\ssup 1}\le t^{q-\e_3}\big\},$$ 
we have the estimate
$$t^{q-\e_2}\leq \Psi_t(Z_t^{\ssup 2})\leq \xi_{\rho_t}^{\ssup 1}+2d\leq \xi_{t^{q+1-\e_1}}^{\ssup 1}+2d\le t^{q-\e_3}+2d,$$
which is impossible for any sufficiently large $t$ since $\e_2<\e_3$. This finishes the proof of~(i).
\smallskip

(ii) Note that 
$$\sum_{z\in\Z^d} u_1^{\ssup t}(t,z)= 
\E_0\Big[\exp\Big\{\int_0^t\xi(X_s)\dd s\Big\}\1\{\tau_{B_{R_t}^{\mathrm{c}}}\le t\}\Big].$$
Again we denote by $J_t$ the number of jumps of the random walk $(X_s \colon s\ge 0)$ 
before time~$t$. As in the proof of~\eqref{est}, using that $J_t\ge R_t$, we obtain 
\begin{equation}\label{max}
\log \sum_{z\in\Z^d} u_1^{\ssup t}(t,z)
\le \log \Big[\sum_{r\ge R_t}e^{t\xi_r^{\ssup 1}-2dt}\frac{(2dt)^r}{r!}\Big]\le 
t\max_{r\ge R_t}\left[\xi_r^{\ssup 1}-\frac{r}{t}\log\frac{r}{2det}\right]+O(t). 
\end{equation}
Denote by $\overline{\rho}_t$ the radius in $\N\cap[R_t,\infty)$ at which the maximum on the right 
hand side is attained, and by $\overline{z}_t$ the maximum point of $\xi$ in the box $B_{\overline \rho_t}$, i.e., the point satisfying 
$|\overline{z}_t|\le \overline{\rho}_t$ and $\xi(\overline{z}_t)=\xi_{\overline{\rho}_t}^{\ssup 1}$. 
\smallskip

Let $f_t=(\log t)^{-\delta}$ and consider the event $\{|Z_t^{\ssup 1}|\ge r_tf_t\}\cap\{\xi(Z_t^{\ssup 1})=\xi_{R_t}^{\ssup 1}\}$.
By Lemma~\ref{psi} and Lemma~\ref{xi} the probability of this converges to one, and so it is sufficient to 
prove the desired estimate on this set. For large~$t$ we have $\overline\rho_t\ge R_t\ge |Z_t^{\ssup 1}|\ge r_tf_t \ge 2dt$.
Supposing for the moment that $|\overline{z}_t|<\overline{\rho}_t$, we obtain, using that
$r\mapsto \frac rt \log \frac{r}{2det}$ is positive and strictly increasing on the interval $(2dt,\infty)$,
\begin{align*}
\xi_{|\overline{z}_t|}^{\ssup 1}-\frac{|\overline z_t|}{t}\log\frac{|\overline z_t|}{2det}
=\xi_{\overline{\rho}_t}^{\ssup 1}-\frac{|\overline z_t|}{t}\log\frac{|\overline z_t|}{2det}
>\xi_{\overline{\rho}_t}^{\ssup 1}-\frac{\overline{\rho_t}}{t}\log\frac{\overline{\rho_t}}{2det},
\end{align*}
which implies $|\overline{z}_t|<R_t$ by definition of $\overline\rho_t$.  Hence
either $|\overline{z}_t|=\overline{\rho}_t$ holds, or $|\overline{z}_t|<R_t$. 
\smallskip
 
In the case $|\overline{z}_t|=\overline{\rho}_t$ we have
\begin{equation}
\label{max1}
\xi_{\overline{\rho}_t}^{\ssup 1}-\frac{\overline{\rho}_t}{t}\log\frac{\overline{\rho}_t}{2det}
=\xi(\overline{z}_t)-\frac{|\overline z_t|}{t}\log\frac{|\overline z_t|}{2det}=\Psi_t(\overline{z}_t)\le \Psi_t(Z_t^{\ssup 2}),
\end{equation}
where the last inequality follows from the fact that $|\overline{z}_t|=\overline{\rho_t}\ge R_t>|Z_t^{\ssup 1}|$ and so 
$\overline{z}_t\neq Z_t^{\ssup 1}$.
\smallskip

In the case $|\overline{z}_t|<R_t$ we use the condition $\xi(Z_t^{\ssup 1})=\xi_{R_t}^{\ssup 1}$ and get
\begin{equation}
\label{max2}
\xi_{\overline{\rho}_t}^{\ssup 1}-\frac{\overline{\rho}_t}{t}\log\frac{\overline{\rho}_t}{2det}
=\xi(\overline z_t)-\frac{\overline{\rho}_t}{t}\log\frac{\overline{\rho}_t}{2det}
\leq \xi_{R_t}^{\ssup 1}-\frac{R_t}t\log\frac{R_t}{2det}
=\xi(Z_t^{\ssup 1})-\frac{R_t}t\log\frac{R_t}{2det}.
\end{equation}

Combining~\eqref{max},~\eqref{max1}, \eqref{max2} we obtain, on the event 
$\{|Z_t^{\ssup 1}|\ge r_tf_t\}\cap\{\xi(Z_t^{\ssup 1})=\xi_{R_t}^{\ssup 1}\}$,
\begin{align*}
\log \sum_{z\in\Z^d} u_1^{\ssup t}(t,z)
\le \max\left\{t\Psi_t(Z_t^{\ssup 2}),t\xi(Z_t^{\ssup 1})-R_t\log\frac{R_t}{2det}\right\}+O(t).
\end{align*}
This completes the proof.
\qed\end{Proof}
\medskip

\begin{Proof}{Proof of Proposition~\ref{lu1u2}.}
Recall the random variables $Y_1\ge Y_2$ from Lemma~\ref{psi2}. Since their joint distribution
is continuous, we have $\P(Y_1=Y_2)=0$. Fix some function $t\mapsto\eta_t$ tending to 0 as $t\to\infty$ 
(to be determined later), then we have
\begin{align}
\label{eta}
\lim_{t\to\infty}\Prob\left(\Psi_t(Z_t^{\ssup 1})-\Psi_t(Z_t^{\ssup 2})\ge {a}_t\eta_t\right)=1.
\end{align}

Fix  $0<\delta<1/4$ and put $f_t=(\log t)^{-\delta}$ and $g_t=(\log t)^{\delta}$. Recall $q=d/(\alpha-d)$ and 
the scale functions $r_t=(t/\log t)^{q+1}$ and $a_t=(t/\log t)^q$. Consider the event
\begin{align*}
\Lambda_t= & \, \Big\{\log \sum_{z\in\Z^d}u_2^{\ssup t}(t,z)\le 
t\Psi_t(Z_t^{\ssup 2})+O(t)\Big\}\\
&\cap\Big\{\log \sum_{z\in\Z^d}u_1^{\ssup t}(t,z)
\le \max\{t\Psi_t(Z_t^{\ssup 2}),t\xi(Z_t^{\ssup 1})-R_t\log\sfrac{R_t}{2det}\}+O(1)\Big\}\\
&\cap\Big\{\log U(t)\ge t\xi(Z_t^{\ssup 1})-|Z_t^{\ssup 1}|\log\xi(Z_t^{\ssup 1})+O(t)\Big\}
\cap\Big\{r_tf_t\le |Z_t^{\ssup 1}|\le r_tg_t\Big\}\\[1mm]
&\cap \Big\{\Psi_t(Z_t^{\ssup 1})\le {a}_t g_t\Big\}
\cap \Big\{\Psi_t(Z_t^{\ssup 1})-\Psi_t(Z_t^{\ssup 2})\ge {a}_t\eta_t\Big\}
\end{align*}
Then $\lim_{t\to\infty}\Prob(\Lambda_t)=1$, according to Lemmas~\ref{bounds2}, \ref{bounds}, \ref{psi} and 
\eqref{psiz} and \eqref{eta}, respectively.\smallskip

On the set $\Lambda_t$ we have the following estimates. \emph{First},
\begin{equation}
\label{u2}
\begin{aligned}
\log\frac{\sum_x u_2^{\ssup t}(t,x)}{U(t)}
&\le t\Psi_t(Z_t^{\ssup 2})-t\xi(Z_t^{\ssup 1})+|Z_t^{\ssup 1}|\,\log\xi(Z_t^{\ssup 1})+O(t)\\
&=-t\left(\Psi_t(Z_t^{\ssup 1})-\Psi_t(Z_t^{\ssup 2})\right)
+|Z_t^{\ssup 1}|\,\log\big[2de\sfrac{t\xi(Z_t^{\ssup 1})}{|Z_t^{\ssup 1}|}\big]+O(t).
\end{aligned}
\end{equation}
For the first of the two terms we get,
$$
-t\left(\Psi_t(Z_t^{\ssup 1})-\Psi_t(Z_t^{\ssup 2})\right)\le -\eta_t t{a}_t =  -\frac{\eta_t t^{q+1}}{(\log t)^{q}},
$$
and for the second,
\begin{equation}\label{u2a}
\begin{aligned}
|Z_t^{\ssup 1}|&\,\log\big[2de\sfrac{t\xi(Z_t^{\ssup 1})}{|Z_t^{\ssup 1}|}\big]+O(t)
\le|Z_t^{\ssup 1}|\,\log\big[\sfrac{t\Psi_t(Z_t^{\ssup 1})}{|Z_t^{\ssup 1}|}
+\log\sfrac{|Z_t^{\ssup 1}|}{2det}\big]+O(|Z_t^{\ssup 1}|)\\
&\le r_tg_t\,\log\big[\sfrac{t{a}_t {g}_t}{r_tf_t}
+\log\sfrac{r_tg_t}{2det}\big]+O(r_tg_t)
\le \sfrac{t^{q+1}\log\log t}{(\log t)^{q+1-\delta}}\,(1+2\delta)\,(1+o(1)).
\end{aligned}
\end{equation}
Now it is clear that we may pick $\eta_t\downarrow0$ such that the right hand side of \eqref{u2} diverges to~$-\infty$.
\smallskip

\emph{Second},
\begin{equation}
\label{u1}
\begin{aligned}
\log\frac{\sum_z u_1^{\ssup t}(t,z)}{U(t)}
&\le \max\left\{t\Psi_t(Z_t^{\ssup 2}),t\xi(Z_t^{\ssup 1})-R_t\log\frac{R_t}{2det}\right\}\\
&-t\xi(Z_t^{\ssup 1})+|Z_t^{\ssup 1}|\log\xi(Z_t^{\ssup 1})+O(t).
\end{aligned}
\end{equation}
We have already shown in~\eqref{u2} that the expression produced 
by the first option in the maximum converges to $-\infty$
on $\Lambda_t$. It remains to show that the same is true for the second option, i.e., for
$$-R_t\,\log\frac{R_t}{2det} + |Z_t^{\ssup 1}|\,\log \xi(Z_t^{\ssup 1}) + O(t)\, .$$
Recalling from \eqref{h} and \eqref{Rtdef} that $R_t=|Z_t|(1+h_t)>|Z_t^{\ssup 1}|$, we obtain 
an upper bound of
$$-h_t|Z_t^{\ssup 1}|\,\log\frac{R_t}{2det}+|Z_t^{\ssup 1}|\,
\log\big[2de\sfrac{t\xi(Z_t^{\ssup 1})}{|Z_t^{\ssup 1}|}\big]+O(t).$$
The first term is  estimated by
$$-h_t|Z_t^{\ssup 1}|\log\frac{R_t}{2det}\le -h_tr_tf_t\log \frac{r_tf_t}{2det}=-\frac{qh_t t^{q+1}}{(\log t)^{q+\delta}}(1+o(1)),$$
while the second is estimated in \eqref{u2a}. One observes that the sum of these two upper bounds diverges to $-\infty$, 
provided that $\frac{(\log t)^{1-2\delta}h_t}{\log\log t}\to \infty$.
This follows from our assumption  $\delta<1/4$ and the definition of $h_t$ in~\eqref{h}. 
Hence the right hand side of~\eqref{u1} goes to $-\infty$. 
\qed\end{Proof}

\subsection{Proof of Proposition~\ref{lu3}: Estimating $\boldsymbol{u_3^{\ssup t}}$.}\label{sec-u3esti}
\\[-2mm]

In this section we prove Proposition~\ref{lu3}, i.e., we show that the total mass of $u_3$ is concentrated on $Z_t$. 
Denote by $\lt$ and $\vt$ the principal eigenvalue and the corresponding 
positive eigenfunction of $\Delta+\xi$ in the box $B_{R_t}$ with zero boundary condition. 
We assume that $\vt$ is normalised to $\vt(Z_t)=1$ and not, as more common, in the~$\ell^2$-sense. 
Then we have the following probabilistic representation of $\vt$,
\begin{align}
\label{ef}
\vt(z)=\E_z\left[\exp\Big\{\int_0^{\tau_{Z_t}}\left(\xi(X_s)-\lt\right) \, \dd s\Big\}
\1\{\tau_{Z_t}<\tau_{B_{R_t}^{\mathrm{c}}}\}\right],
\qquad \mbox{ for } z\in B_{R_t}^{\mathrm{c}}.
\end{align}
\begin{lemma} 
\label{lem-lu3}\ \\[-3mm]
\begin{enumerate}
\item[(i)] For any $t>0$ and all $(\theta,z)\in (0,\infty)\times\Z^d$, we have
\begin{align*}
u_3^{\ssup t}(\theta,z)\le u_3^{\ssup t}(\theta,Z_t) \, ||\vt||_2^2 \, \vt(z).
\end{align*}
\item[(ii)] The eigenfunction $\vt$ is localised around $Z_t$ so that
\begin{align*}
||\vt||_2^2 \sum_{z\neq Z_t}\vt(z)\to 0\quad\text{ in probability}.
\end{align*}
\end{enumerate}
%in probability as $t\to\infty$.
\end{lemma}

\begin{Proof}{Proof.}
(i) The first estimate is a special case of~\cite[Th.~4.1]{GKM06} (with $\Gamma=\{Z_t\}$
in their notation) but we repeat the proof here for the sake of completeness. By time reversal, 
\begin{align*}
u_3^{\ssup t}(\theta,z)&=\E_z\left[\exp\Big\{\int_0^{\theta}\xi(X_{s})\, \dd s\Big\}\, \1\{X_{\theta}=0\}
\, \1\{\tau_{B_{R_t}^{\mathrm{c}}}>\theta\}\1\{\tau_{Z_t}\le \theta\}\right].
\end{align*}
We obtain a lower bound for $u_3^{\ssup t}(\theta,Z_t)$
by requiring that the random walk is in $Z_t$ at time $u\in (0,\theta)$. Using the Markov property 
at time $u$, we obtain 
\begin{align}
u_3^{\ssup t}(\theta,Z_t)&\ge 
\E_{Z_t}\left[\exp\Big\{\int_0^u
\xi(X_{s})\, \dd s\Big\}\,\1\{Z_t=X_u\}\,
\1\{\tau_{B_{R_t}^{\mathrm{c}}}>u\}\right]\notag\\
&\qquad \times
\E_{Z_t}\left[
\phantom{\int}\right.\!\!\!\!\!\!\!
\exp\Big\{\int_0^{\theta-u}\xi(X_{s})\dd s\Big\}\,\1\{X_{\theta-u}=0\}\,
\1\{\tau_{B_{R_t}^{\mathrm{c}}}>\theta-u\}
\!\!\!\!\!\!\!\!\left.\phantom{\int}
\right].
\label{st1}
\end{align}  
Using an eigenvalue expansion for the parabolic problem in $B_{R_t}$ represented by the 
first factor in the formula above, we obtain the bound
\begin{align*}
\E_{Z_t}\left[\exp\Big\{\int_0^u
\xi(X_{s})\, \dd s\Big\} \, \1\{Z_t=X_u\}
\,\1\{\tau_{B_{R_t}^{\mathrm{c}}}>u\}\right]
\ge \ee^{\lt u}\, \frac{\vt(Z_t)^2}{||\vt||_2^2}
=\ee^{\lt u} \, ||\vt||_2^{-2}.
\end{align*}
Substituting this into~\eqref{st1}, for $0<u<\theta$,
\begin{align}
\label{st2}
\E_{Z_t}\left[
\phantom{\int}\right.\!\!\!\!\!\!\!
\exp\Big\{\int_0^{\theta-u}\xi(X_{s})\, \dd s\Big\}\,\1\{X_{\theta-u}=0\}
\,\1\{\tau_{B_{R_t}^{\mathrm{c}}}>\theta-u\}
\!\!\!\!\!\!\!\!\left.\phantom{\int}
\right]
&\le\ee^{-\lt u} \, ||\vt||_2^2 \, u_3^{\ssup t}(\theta,Z_t).
\end{align} 
Since the claimed estimate
%~\eqref{u3_1} 
is obvious for $z=Z_t$ due to the norming of $\vt$, we may assume that
$z\in B_{R_t}\setminus\{ Z_t\}$. 
Using the strong Markov property at time $\tau_{Z_t}$ 
and~\eqref{st2} with $u=\tau_{Z_t}$
we obtain
\begin{align*}
u_3^{\ssup t}(\theta,z)&
=\E_z\left[\exp\Big\{\int_0^{\tau_{Z_t}}
\xi(X_{s})\,\dd s\Big\} \, 
\1\{\tau_{B_{R_t}^{\mathrm{c}}}>\tau_{Z_t}\} \, \1\{\tau_{Z_t}\le \theta\}\right.\\
&\phantom{aaaaa}\times
\E_{Z_t}\left[\phantom{\int}\right.\!\!\!\!\!\!\exp\Big\{\int_0^{\theta-u}\xi(X_{s}) \, \dd s\Big\}
\, \1\{X_{\theta-u}=0\} \, \1\{\tau_{B_{R_t}^{\mathrm{c}}}>\theta-u\}\!\!\!\!\!\!\left.\phantom{\int}\right]_{u=\tau_{Z_t}}
\left.\phantom{\int_0^{\tau_{\overline{Z}_t^{\ssup 1}}}}\!\!\!\!\!\!\!\!\!\!\!\!\!\!\!\!\right]\\
&\le
u_3^{\ssup t}(\theta,Z_t) \, ||\vt||_2^2 \, \E_z\left[\exp\Big\{\int_0^{\tau_{Z_t}}
(\xi(X_{s})-\lt)\, \dd s\Big\} \, 
\1\{\tau_{Z_t}<\tau_{B_{R_t}^{\mathrm{c}}}\}\right]\\
&=u_3^{\ssup t}(\theta,Z_t) \, ||\vt||_2^2 \, \vt(z).
\end{align*}

(ii) To prove the localisation of $\vt$ around $Z_t$, first note that
by the Rayleigh-Ritz formula
\begin{equation}
\label{rr}
\begin{aligned}
\lt
&=\sup\left\{\langle (\Delta+\xi)f,f\rangle\colon f\in \ell^2(\Z^d),\mathrm{supp}(f)\subset B_{R_t},||f||_2=1\right\}\\
&\ge \sup\left\{\langle (\Delta+\xi)\delta_z,\delta_z\rangle\colon z\in B_{R_t}\right\}
=\sup\left\{\xi(z)-2d\colon z\in B_{R_t}\right\}\\
&=\xi_{R_t}^{\ssup 1}-2d.
\end{aligned}
\end{equation}
Recall the diverging  function $b\colon(0,\infty)\to(0,\infty)$ from Lemma~\ref{xi}(iv) and consider the event
$$\Lambda_t=\left\{\xi_{R_t}^{\ssup 1}-\xi_{R_t}^{\ssup 2}\ge{b}_t\right\}
\cap \left\{\xi(Z_t)=\xi_{R_t}^{\ssup 1}\right\}.$$
By Lemma~\ref{xi}(ii) and (iv), its probability converges to $1$.
It follows from~\eqref{rr} that, on $\Lambda_t$,
$$\xi_{R_t}^{\ssup 2}-\lt\le \xi_{R_t}^{\ssup 1}-{b}_t-\lt \le 2d-{b}_t.$$
Since the paths of 
the random walk $(X_s \colon s\ge 0)$ in~\eqref{ef} do not leave $B_{R_t}$ and 
avoid the point~$Z_t$ where the maximum $\xi_{R_t}^{\ssup 1}$ is achieved, we can estimate the integrand 
in terms of the second-largest value of $\xi$ in $B_{R_t}$. Hence, we obtain 
\begin{align}
\vt(z)
&\le \E_z\big[\exp\big\{\tau_{Z_t}\left(\xi_{R_t}^{\ssup 2}-\lt\right)\big\}\big]
\le \E_z\big[\exp\big\{\tau_{Z_t}\left(2d-{b}_t\right)\big\}\big].
\label{estv}
\end{align}
Under $\P_z$ the random variable $\tau_{Z_t}$ is stochastically bounded from below by a sum of $|z-Z_t|$ independent exponentially 
distributed random times with parameter $2d$. If $\tau$ denotes such a random time, we therefore have
$$
\vt(z)\leq \E_z\Big[\exp\big\{\tau_{Z_t}(2d-{b}_t)\big\} \Big]
\leq \Big(\E\big[e^{-[b_t-2d]\tau}\big]\Big)^{|z-Z_t|}
=\Big(\frac{2d}{b_t}\Big)^{|z-Z_t|}.
$$
From this, it is easy to see that the assertion holds.
\qed\end{Proof}
\smallskip

\begin{Proof}{Proof of Proposition~\ref{lu3}.} 
Lemma~\ref{lem-lu3}(i) yields that
\begin{align*}
\frac{\sum_{z\neq Z_t} u_3^{\ssup t}(t,z)}{U(t)}\le \frac{\sum_{z\neq Z_t}
u_3^{\ssup t}(t,z)}{u_3^{\ssup t}(t,Z_t)}
%=\sum\frac{u_3^{\ssup t}(t,z)}{u_3^{\ssup t}(t,Z_t)}
\le ||\vt||_2^2\sum_{z\neq Z_t}\vt(z),
\end{align*}
and the right hand side vanishes as $t\to\infty$
in probability, by Lemma~\ref{lem-lu3}(ii).
\qed\end{Proof}

\bigskip\noindent{\bf Acknowledgements:} We thank Remco van der Hofstad for helpful
discussions. The first author is supported by the Forschergruppe 718 of the German Science Foundation (DFG). 
The second author is supported by an Advanced Research Fellowship, and both the second and third author by grant EP/C500229/1 of the Engineering and Physical Sciences Research Council (EPSRC).

\vspace{0.4cm}

{\footnotesize
\begin{tabular}{ll}
{\bf Wolfgang K\"onig} & {\bf Peter M\"orters} \\
Universit\"at Leipzig & University of Bath\\
Fakult\"at f\"ur Mathematik und Informatik & Department of Mathematical Sciences\\
Mathematisches Institut & Claverton Down\\
Augustusplatz 10/11 &  Bath \\
04109 Leipzig & BA2 7AY\\
Germany & United Kingdom\\
{\tt koenig@math.uni-leipzig.de} & {\tt maspm@bath.ac.uk}\\
& \\[0.1cm]
{\bf Nadia Sidorova} & \\
University of Bath & \\
Department of Mathematical Sciences & \\
Claverton Down & \\
Bath & \\
BA2 7AY & \\
United Kingdom & \\
{\tt n.sidorova@maths.bath.ac.uk} & \\
\end{tabular}}

\end{document}